\documentclass[10pt,twoside,twocolumn]{article}

\usepackage[margin=0.8in,columnsep=0.3in]{geometry}
\usepackage[utf8]{inputenc}
\usepackage[T2A]{fontenc}
\usepackage[russian,english]{babel}
\usepackage{amsmath,amssymb,amsthm,mathtools}
\usepackage{algorithm}
\usepackage{algpseudocode}
\usepackage{booktabs}
\usepackage{tabularx,array}
\usepackage{graphicx}
\usepackage{xcolor}
\usepackage[hidelinks]{hyperref}
\usepackage[capitalize,noabbrev]{cleveref}
\usepackage[affil-it]{authblk}

\newcommand{\citep}{\cite}
\newcommand{\citet}{\cite}

\theoremstyle{plain}
\newtheorem{theorem}{Theorem}
\newtheorem{lemma}{Lemma}
\newtheorem{proposition}{Proposition}

\theoremstyle{definition}
\newtheorem{definition}{Definition}
\newtheorem{assumption}{Assumption}
\theoremstyle{remark}
\newtheorem{remark}{Remark}

\crefname{assumption}{Assumption}{Assumptions}
\crefname{algorithm}{Algorithm}{Algorithms}

\newcommand{\R}{\mathbb{R}}
\newcommand{\E}{\mathbb{E}}
\renewcommand{\Pr}{\mathbb{P}}
\newcommand{\B}{\mathbb{B}}
\renewcommand{\S}{\mathbb{S}}
\newcommand{\cA}{\mathcal{A}}
\newcommand{\cF}{\mathcal{F}}
\newcommand{\cH}{\mathcal{H}}
\newcommand{\cL}{\mathcal{L}}
\newcommand{\cN}{\mathcal{N}}
\newcommand{\cO}{\mathcal{O}}
\newcommand{\Fm}{\mathcal{F}_\mu^M}
\newcommand{\Fzero}{\mathcal{F}_0^M}
\DeclareMathOperator{\sign}{sign}
\DeclareMathOperator{\Projp}{Proj}
\DeclareMathOperator{\Unif}{Unif}

\DeclareMathOperator{\KL}{KL}
\DeclareMathOperator{\sech}{sech}
\DeclareMathOperator{\TV}{TV}
\newcommand{\ind}{\mathbf{1}}
\newcommand{\ip}[2]{\langle #1,#2\rangle}

\title{Nonsmooth Convex and Strongly Convex Optimization \\from Noisy Comparisons}
\author[1]{Evgeniy Kovalev}
\author[1, 2]{Fedor Stonyakin}
\author[1,2]{Alexander Gasnikov}
\affil[1]{Innopolis University, Russia}
\affil[2]{Moscow Institute of Physics and Technology, Russia}
\date{}
\date{}

\begin{document}

\maketitle

\begin{abstract}
We study nonsmooth convex and strongly convex optimization problems using only noisy pairwise comparisons: each query specifies two points and returns a single binary response indicating which point has the smaller function value, subject to noise. We distinguish two observation models: \emph{Stochastic Order Oracles}, which compare two evaluations of the same random loss, and \emph{Preference Oracles}, whose conditional response mean is a transfer function of the objective difference. For Stochastic Order Oracles, we prove that no comparison-based algorithm can guarantee vanishing expected suboptimality in general, even when sample losses are convex and Lipschitz and the mean objective is strongly convex. Under a normalized probit calibration, we recover convergence and establish matching upper and lower rates in the query budget for convex and strongly convex objectives. For Preference Oracles with a linear response model, the binary response $O \in \{\pm1\}$ satisfies $\E[O\mid y,z]=\lambda(f(y)-f(z))$, with $\lambda>0$. In dimension $d$, we obtain an upper bound of $\widetilde O(d^{3/2}/(\lambda\sqrt N))$ and establish, to our knowledge, the first minimax lower bound of $\Omega(d/(\lambda\sqrt N))$, matching the optimal budget $N$ exponent, for both convex and strongly convex nonsmooth objectives. Beyond exact linearity, we analyze a simple projected comparison method for nonlinear transfer functions and establish convergence rates that depend on how accurately these functions can be approximated by a linear model near zero. Our analysis requires neither oddness nor twice differentiability of the transfer function.
 
\end{abstract}

\section{Introduction}
\label{sec:intro}

In many learning systems the most reliable feedback is \emph{relative}. Human raters
compare two responses more consistently than they score one
\citep{thurstone1927law,bradley1952rank}; reward models for language models are trained
from pairwise preferences \citep{christiano2017deep,ouyang2022training}; and automatic
judges are typically queried with pairs \citep{zheng2023judging}. The same structure
appears in preference-based Bayesian optimization \citep{gonzalez2017preferential},
interactive information retrieval \citep{yue2009interactively}, and robot reward learning
\citep{sadigh2017active}. This motivates the optimization problem
\begin{equation}\label{eq:problem}
  \min_{x\in Q} f(x),\qquad Q\subset\R^d\text{ compact and convex},
\end{equation}
in which an algorithm may query pairs $(x,y)\in Q\times Q$ and observes only a noisy bit
indicating which point appears better. We ask how many comparisons are needed to reach
accuracy $\varepsilon$ when $f$ is convex and Lipschitz but not smooth.
Non-smoothness is the rule rather than the exception in this setting: utilities with
clipping or thresholds, worst-case and $\ell_1$-type regularizers are all non-differentiable.

\paragraph{Two observation models.}
What a comparison measures depends on the source of noise. In a
\emph{Stochastic Order Oracle}, firstly introduced in \cite[Section~5]{lobanov2024acceleration} and studied in \cite{smirnov2024asymptotic}, the
problem is
\begin{equation*}
    \min_{x\in Q}f(x),\quad f(x)=\mathbb{E}_\xi[f(x,\xi)],\quad Q\subset \mathbb{R}^d,
\end{equation*}
and the oracle reveals
$\sign(f(x,\xi)-f(y,\xi))$ for a fresh realization $\xi$ shared by both points. 
In a
\emph{Preference Oracle} \citep{saha2025dueling,lu2026online}, the problem is
\begin{equation*}
    \min_{x\in Q}f(x),\quad Q\subset \mathbb{R}^d,
\end{equation*}
where a response of the oracle is a bit
with $\E[O\mid x,y]=\psi(f(x)-f(y))$ for a transfer function $\psi$, as in the
Thurstone \cite{thurstone1927law} and Bradley--Terry \cite{bradley1952rank} models. Despite that two models look alike, they behave very
differently. It happens due to the fact that in the first scenario, the function is stochastic itself and the oracle is precise, meanwhile in the second setting, the function is determined, but the noise is in the oracle. More differences will be seen in the following sections.

\paragraph{Our contribution.}
\begin{itemize}
\item \textbf{Divergence for Stochastic Order Oracle} (\cref{thm:nonident}). 
In the first setting, where the algorithm has access only to a Stochastic Order Oracle, 
we show that no comparison-based algorithm can guarantee vanishing expected error, even 
when every sample loss is convex and $M$-Lipschitz and the mean is strongly convex. 
The obstruction is that comparisons are invariant to positive random rescalings of the 
sample loss, yet such rescalings shift the minimizer of the mean.
\item \textbf{Convergence in particular cases} (\cref{thm:probit}). Nevertheless, when the losses 
exhibit additive Gaussian linear noise (which is the random-coefficient probit model of discrete 
choice), we show that Algorithm \ref{alg:probit} attains the rate
\[
\Theta\!\left(\min\left\{\,d/\sqrt N,\; d^2/(\mu N)\right\}\right),
\]
where the second term applies only under strong convexity. For explicit constants, refer to \cref{thm:probit}.
\item \textbf{First lower and upper bounds for the Preference Oracle} 
(\cref{thm:lin-upper,thm:lin-lower}). In the second setting, where the algorithm has 
access only to a Preference Oracle with linear $\psi(t)=\lambda t$ for some $\lambda>0$, we 
prove an upper bound of $\widetilde O\!\left(d^{3/2}/(\lambda\sqrt N)\right)$ and a 
lower bound of $\Omega\!\left(d/(\lambda\sqrt N)\right)$ for convex and strongly convex 
non-smooth $f$, matching in the exponent of $N$. To our knowledge, this is the first 
minimax lower bound for comparison feedback in the non-smooth setting, and it shows 
that strong convexity cannot improve upon the $N^{-1/2}$ rate.
\item \textbf{First lower and upper bounds for general transfer functions} 
(\cref{thm:ardg-convex,thm:ardg-sc}). We propose \cref{alg:ardg}, a modification of a 
classic projected comparison method \cite{saha2025dueling}, and prove that it converges 
at rates governed only by the local linearization exponent $q$ of $\psi$. We further show that these rates are not improvable among methods 
with the same radius budget (\cref{thm:restricted}).

\item \textbf{Experiments} (\cref{sec:experiments}). We confirm the impossibility result and
the rates, including for transfer functions outside the scope of prior analyses.
\end{itemize}

\paragraph{Related work.}
\emph{Dueling bandits} have been studied for comparisons over finite or structured
action sets in \cite{yue2012karmed,bengs2021survey}. For continuous convex problems,
dueling bandit gradient descent was proposed for this setting in \cite{yue2009interactively},
while matching upper and lower bounds were proven for strongly convex smooth objectives
under power-law transfer functions - in \cite{jamieson2012query}. Regret bounds for strongly
convex smooth objectives were later established in \cite{kumagai2017regret}, and the
noiseless sign oracle was treated by authors of \citep{saha2021dueling}. Then, by the same researchers, a setting for smooth $f$ and general odd and differentiable transfer functions was examined in \cite{saha2025dueling}.
Adversarial online optimization with dueling feedback was studied in \cite{lu2026online}, assuming an odd $\psi$ with bounded $\psi''$. In contrast, we handle non-smooth $f$, give the first lower bound without smoothness, and
allow non-odd, non-sign-preserving, and non-differentiable $\psi$ (\cref{tab:comparison}). Our analysis builds
on two-point gradient estimation from bandit convex optimization
\citep{flaxman2005online,agarwal2010optimal,duchi2015optimal,shamir2017optimal} and on
ball smoothing \citep{nemirovski1983problem,gasnikov2022power}. Our linear-link bounds
inherit the $\sqrt d$ gap between the upper bound of prior work \citep{fokkema2024online}
and the lower bound known from \citep{shamir2013complexity}. Deterministic order oracles
have been studied in \cite{bergou2020stochastic,lobanov2024acceleration}, and ranking
oracles for human feedback were used  in \citep{tang2024zeroth}. Rank-based
evolution strategies \citep{hansen2001completely} are invariant to monotone
transformations of the loss, which is the same invariance that underlies
\cref{thm:nonident}.

\section{Setup and Main Definitions}
\label{sec:setup}
Before presenting our results, we introduce the notation and main assumptions we use in our work.

\begin{definition}[Function classes]\label{def:classes}
For $\mu\ge0$ and $M>0$, $\Fm(Q)$ is the class of functions $f:Q\to\R$ that are
$M$-Lipschitz and $\mu$-strongly convex on $Q$; $\Fzero(Q)$ is the convex case. When obvious, we omit dependency of $Q$ for simplicity.
\end{definition}

Throughout the paper, $\|\cdot\|$ denotes the Euclidean norm, and $\B^d(c,r)$ and 
$\S^{d}(c,r)$ denote the Euclidean ball and sphere centered at $c$ with radius $r$. 
Unless stated otherwise, $Q$ is convex and compact with diameter $D<\infty$. We denote 
by $\Projp_Q$ the Euclidean projection onto $Q$. When obvious, we omit the dimension $d$ 
to simplify notation.

\begin{definition}\label{def:alg-class}
An algorithm $A\in\cA_N$ consists of a probability space $(\mathcal U,\nu)$ (its internal
randomness) and measurable maps $\phi_1,\dots,\phi_{N+1}$. Run on an oracle $\mathcal O^f$ for
an objective $f$ (we write $\mathcal O$ when $f$ is clear), it draws a seed $U\sim\nu$
independent of the oracle and, for $t=1,\dots,N$, queries the pair
\[
  (x_t,y_t):=\phi_t(U,\mathcal{O}_1,\dots\mathcal{O}_{t-1})\in Q\times Q
\]
and receives the answer $\mathcal{O}_t:=\mathcal O(x_t,y_t)$. Its output is
$\hat x:=\phi_{N+1}(U,\mathcal{O}_1,\dots,\mathcal{O}_N)\in Q$. An algorithm that makes fewer than $N$ calls is
included by padding it with calls whose answers it ignores.
\end{definition}

\begin{definition}[Minimax risk]
For an oracle model $\cO$ and a class $\cF$, we will define the minimax risk as
\[
  \cL_N(\cF,\cO)\;:=\;\inf_{A\in\cA_N}\ \sup_{f\in\cF}\ \E\big[f(\hat x_N)-f^*\big],
  \quad f^*:=\min_{x\in Q} f.
\]
That is, it is the best error guarantee achievable over $\cF$ by any algorithm in
$\cA_N$ that makes at most $N$ oracle calls.

\end{definition}
The central technique underlying our proposed algorithms is $\ell_2$-ball smoothing, 
which dates back to \citet{nemirovski1983problem} and \citet{gasnikov2022power} and was 
studied further in \citet{kornilov2023accelerated,gasnikov2022power}. More precisely, for $h>0$, let $f_h(x):=\E_{u\sim\Unif(\B^d)}f(x+hu)$.

The power of this technique is captured by the following lemma, established in \cite{nemirovski1983problem,flaxman2005online}, which provides access to 
a gradient of the objective through its smoothed proxy. For completeness, we also provide a proof of that lemma in \cref{proof:smoothing}.
\begin{lemma}\label{lem:smoothing}
If $f$ is convex and $M$-Lipschitz on the $h$-neighborhood of $Q$, then on $Q$:
\begin{enumerate}
  \item $f_h$ is convex, and $\mu$-strongly convex if $f$ is;
  \item $0 \le f_h - f \le Mh$;
  \item $f_h$ is differentiable with
  \[
    \nabla f_h(x) = \frac{d}{2h}\,\E_{e\sim\Unif(\S^{d-1})}
    \big[(f(x+he)-f(x-he))\,e\big].
  \]
\end{enumerate}
\end{lemma}
\section{Stochastic Order Oracles}
\label{sec:order}
In this scenario, we consider the optimization problem 
\begin{equation}
    \min_{x\in Q}f(x),\quad f(x)=\mathbb{E}_\xi[f(x,\xi)],\quad Q\subset \mathbb{R}^d.
\end{equation}

\begin{definition}[Stochastic Order Oracle]\label{def:order}
At step $t$, the algorithm queries a pair $(X_t,Y_t)\in Q\times Q$ and observes
\[
  O_t = \sign\big(f(X_t,\xi_t)-f(Y_t,\xi_t)\big),
\]
where $\xi_t$ is a fresh sample shared by both points, and ties are broken by an 
independent fair coin.
\end{definition}

This oracle models settings in which two candidates are evaluated under the same random 
conditions. For example, a randomly sampled individual compares two products according 
to their own taste, as motivated by \citet{lobanov2024acceleration}; two simulator 
configurations are run with common random numbers; or two candidates are ranked on the 
same minibatch, as in rank-based evolution strategies \citep{hansen2001completely}.
\begin{theorem}\label{thm:nonident}
Let $Q=\B^d(R)$. Then, for any $N\ge 1$, $M>0$, and $0\le\mu\le M/(8R)$,
\[
\cL_N(\Fm,O_t)\ge\frac{5MR}{32}-\frac{\mu R^2}{2}\;\ge\;\frac{3MR}{32}.
\]
\end{theorem}

The counter-example is built by defining two instances that are nearly indistinguishable. Concretely, let
\(
f_\sigma(x,\xi)=\Big(1+\tfrac34\sigma\xi\Big)\,q(x,\xi),
\)
where $\xi$ takes the values $-1$ and $+1$ with equal probability, and
\(
  q(x,\xi)=\tfrac\mu2\|x\|^2+\tfrac M4\xi x_1+\tfrac M{32}|x_1|.
\)
The positive prefactor is the same for both points in any comparison, so it does not affect the sign of $f_\sigma(X,\xi)-f_\sigma(Y,\xi)$. At the same time, this prefactor is 
correlated with $\xi$ and therefore biases its mean. The proof is given in \ref{app:nonident}.

Therefore, one can see that any guarantee for this oracle requires an assumption that links the law of the comparison to $f$. We now
make two such links explicit.

\subsection{Specific Case I. Probit Comparisons}

Despite the general impossibility, we now show that convergence can be established under 
some specific restrictions on the noise. Suppose the sample loss is perturbed by additive 
Gaussian linear noise,
\begin{equation}\label{eq:probit-model}
f(X,\xi)=f(X)+\langle\xi,X\rangle+g(\xi),\qquad \xi\sim\cN(0,\sigma^2I_d),
\end{equation}
with arbitrary $g$. Then,
\begin{equation}\label{eq:speccase1}
    \Pr(O=1\mid X,Y)=\Phi\big((f(X)-f(Y))/(\sigma\|X-Y\|)\big),
\end{equation} where $\Phi$ is the standard
normal distribution function. In discrete choice, this is the random-coefficient probit
model: individuals share the systematic utility $f$ but have heterogeneous Gaussian linear
tastes \citep{hausman1978conditional,train2009discrete}. In stochastic optimization, it is
the classical additive gradient-noise model $\nabla f(X,\xi)=\nabla f(X)+\xi$ whose asymptotic behaviour was analyzed in \cite{lobanov2024acceleration,smirnov2024asymptotic}.

\begin{theorem}\label{thm:probit}
Let $B(c,r)\subseteq Q\subseteq B(c,R)$, $f\in\Fm$ with $0\le\mu\le M/(8R)$, and assume
\eqref{eq:speccase1} for all distinct $X, Y\in Q$. Put
\[
  \mathcal T_\mu(N,r,s):=\min\Big\{Mr,\ \frac{dsr}{\sqrt N},\ \frac{d^2s^2}{\mu N}\Big\},
\]
omitting the last term if $\mu=0$.

\textup{(a)} Then, for any $N\ge1$, \cref{alg:probit} produces $\hat x_N$, such that \[
  \E[f(\hat x_N)]-f^*\le C(M/\sigma)\,\mathcal T_\mu(N,R,\sigma),
\]
where constant $C(M/\sigma)<\infty$ depends only on $M/\sigma$.

\textup{(b)} Conversely, for every $N\ge1$ and $\mu\ge0$,
\[
  \cL_N(\Fm,O_t)\ge\frac1{100}\,\mathcal T_\mu(N,r,\sigma).
\]
\end{theorem}

The proof and definition of the constant $C(\cdot)$ are given in \cref{proof:thm_probit,proof:thm-probit-b}.

Two remarks are in order. First, the rate is a factor $\sqrt d$ slower than for a stochastic
first-order oracle with the same noise \citep{nemirovski1983problem}, since each query
reveals a single directional difference. Second, the normalization by $\|X-Y\|$ makes a
comparison equally informative at every scale, so the probing radius can be taken small
without a variance penalty, and the bias--variance trade-off of smoothing disappears.

\subsection{Specific Case II: Linear Comparisons}

If instead $\E[O\mid X,Y]=\lambda\big(f(X)-f(Y)\big)$, for instance because the
difference of the sample losses equals $f(x)-f(y)$ plus noise that is uniform on
$[-A,A]$ with $A\ge2MR$, then the order oracle is a Preference Oracle with linear
transfer function, and all results of \cref{sec:linear} apply verbatim. This
calibration leads to the second model.

\section{Preference Oracles}
\label{sec:preference}

In this setting, we consider the optimization problem
\begin{equation}\label{saha_problem}
    \min _{x\in Q}f(x), \qquad Q\subset \mathbb{R}^d.
\end{equation}
\begin{definition}[Preference Oracle]\label{saha_oracle_nonlinear}
At step $t$, the algorithm queries a pair $(x_t,y_t)\in Q\times Q$ and observes
\[
  o_t\in\{\pm1\},\qquad
  \E[o_t\mid x_t,y_t]=\psi\big(f(x_t)-f(y_t)\big),
\]
\end{definition}

for a measurable function
$\psi:\R\to[-1,1]$ unknown to the algorithm. We call $\psi$ the \textit{transfer function}. The value  $\psi(f(x)-f(y))$ is interpreted as the degree to which $x$ should be preferred over $y$ (or vice versa).

The transfer function describes the rater. Thurstone's model with Gaussian utility noise \cite{thurstone1927law} gives
$\psi(t)=2\Phi(t/(\sqrt2\sigma))-1$, and Bradley--Terry gives $\psi(t)=\tanh(t/(2s))$ \cite{bradley1952rank} for some constants $\sigma, s>0$. The latter is the model used to
train reward models from human feedback \cite{christiano2017deep,ouyang2022training}.
Relative judgements (where a rater only reports which of two items is preferred) are easier to elicit and less biased than absolute scores (where a rater assigns one a number on a fixed scale)
\cite{saha2025dueling}, and they underlie language-model alignment, robotics, and
recommender systems \cite{lu2026online}. Since $\psi$ is unknown in practice,
guarantees should depend on it only through local, checkable quantities. \cref{tab:comparison} compares our results  with prior work.

\begin{table*}[t]
\caption{Comparison with prior work on convex optimization from preference
feedback. Rates are shown in terms of $d$ and $N$ only, with other factors omitted. $\widetilde O$ hides logarithmic factors, $q$ is the exponent of \cref{as:q}, and
$\rho(T):=\sup_{|t|\le T}|\psi(t)-\lambda t|$.}
\label{tab:comparison}
{\small
\setlength{\tabcolsep}{4pt}
\renewcommand{\tabularxcolumn}[1]{m{#1}}
\begin{tabularx}{\textwidth}{@{}>{\raggedright\arraybackslash}m{2.4cm}>{\raggedright\arraybackslash}X>{\centering\arraybackslash}m{1.5cm}>{\centering\arraybackslash}m{2.5cm}>{\centering\arraybackslash}m{2.5cm}>{\centering\arraybackslash}m{3.1cm}@{}}
\toprule
Work & Transfer function $\psi$ & Objective & Convex & Strongly convex & Lower bound\\
\midrule
\citet{saha2025dueling}
 & linear
 & smooth
 & $O(d^{2/3}N^{-1/3})$\,$^{\dagger}$
 & $O(d^{3/2}N^{-1/2})$
 & ---\\
 \addlinespace
\citet{saha2025dueling}
 & differentiable, odd, sign-preserving, $\psi'(t)\ge c\,p|t|^{p-1}$, $p\ge1$
 & smooth
 & $O(d^{\frac{2p+1}{4p}}N^{-\frac1{4p}})$\,$^{\dagger}$
 & $O(d^{\frac{2p+1}{2p}}N^{-\frac1{2p}})$
 & ---\\
 \citet{jamieson2012query}
 & Same$\,^{\star}$
 & smooth
 & ---
 & $\widetilde O(d^{1+\frac1{2p}}N^{-\frac1{2p}})$
 & $\Omega(d^{\frac1{2p}}N^{-\frac1{2p}})$\\
\addlinespace
\addlinespace
\citet{lu2026online}
 & odd, $\psi'(0)\neq0$, $|\psi'''|$ bounded
 & smooth
 & $\widetilde O(d\,N^{-1/3})$
 & $\widetilde O(d^{3/2}\,N^{-1/2})$
 & ---\\
\addlinespace
 \citet{lu2026online}
 & odd, $\psi'(0)\neq0$, $|\psi''|$ bounded
 & non-smooth
 & $O(d^{1/2}N^{-1/4})$
 & $\widetilde O(d^{2/3}N^{-1/3})$
 & ---\\
\midrule
\textbf{Ours}, \cref{thm:lin-upper,thm:lin-lower}
 & linear
 & non-smooth
 & $\widetilde O(d^{3/2}N^{-1/2})$
 & $\widetilde O(d^{3/2}N^{-1/2})$
& $\Omega(d\,N^{-1/2})$ \shortstack{convex and \\ strongly convex}\\
\addlinespace
\textbf{Ours}, \cref{thm:ardg-convex,thm:ardg-sc}
 & $\rho(T)\le LT^{1+q}$, $q\ge1$
 & non-smooth
 & $O(d^{3/4}N^{-1/4})$
 & $O(d^{4/3}N^{-1/3})$
 & ---\\
\addlinespace
\textbf{Ours}, \cref{thm:ardg-convex,thm:ardg-sc}
 & $\rho(T)\le LT^{1+q}$, $\frac12\le q<1$
 & non-smooth
 & $O(d^{\frac{1+2q}{2(1+q)}}N^{-\frac{q}{2(1+q)}})$
 & $O(d^{4/3}N^{-1/3})$
 & ---\\
\addlinespace
\textbf{Ours}, \cref{thm:ardg-convex,thm:ardg-sc}
 & $\rho(T)\le LT^{1+q}$, $0<q<\frac12$
 & non-smooth
 & $O(d^{\frac{1+2q}{2(1+q)}}N^{-\frac{q}{2(1+q)}})$
 & $O(d^{\frac{1+2q}{1+q}}N^{-\frac{q}{1+q}})$
 & ---\\
\addlinespace

\bottomrule
\end{tabularx}
}
{\footnotesize
$^{\dagger}$Guarantee on $\min_{t\le N}\E f(w_t)$, i.e.\ for the best iterate, which is not
identified by the algorithm.
$^{\star}$\citet{jamieson2012query} assumes no transfer function: with $t=f(x)-f(y)$, they only
require $\Pr\big(o=\mathrm{sign}(t)\big)\ge\frac12+\min\{\delta_0,\mu|t|^{p}\}$. For
$\psi(t)=c\,\mathrm{sign}(t)|t|^{p}$, this is covered by the model of \cite{jamieson2012query}.
}
\end{table*}

\subsection{Linear Transfer Function}
\label{sec:linear}
The linear link is a special case that arises in many practical situations, since it is equivalent to a bandit feedback optimization setting \cite{saha2025dueling,flaxman2005online}.

\begin{assumption}[Linear link]\label{as:linear-link}
There is $\lambda>0$ with $2\lambda MR\le1$ such that
\begin{equation}\label{eq:linear-link}
  \psi(t)=\lambda t\qquad\text{for }|t|\le2MR,
\end{equation}
where $Q\subseteq \B^d(c,R)$.
\end{assumption}

The constraint $2\lambda MR\le1$ only states that $\lambda(f(y)-f(z))$ is a valid mean
of a $\pm1$ variable on the attainable range of gaps. 

\emph{One-bit two-point feedback.} If an evaluator produces a bounded score difference
$S\in[-1,1]$ with $\E[S\mid x,y]=\lambda(f(x)-f(y))$, which is two-point bandit feedback
\cite{agarwal2010optimal,shamir2017optimal}, and transmits only
$o=+1$ with probability $(1+S)/2$ (dithered or stochastic rounding,
\cite{gray1993dithered,alistarh2017qsgd}), then $o$ follows~\eqref{eq:linear-link}.
The linear-link model is therefore precisely zeroth-order optimization compressed to one
bit, and our bounds quantify the price of this compression. This is the sense in which
\citet{saha2025dueling} call it equivalent to bandit convex optimization
\citep{flaxman2005online}.

\emph{Bounded additive noise.} If $o=\sign(f(x)-f(y)+\xi)$ with $\xi\sim\Unif[-A,A]$ and
$A\ge2MR$, then $\E[o\mid x,y]=(f(x)-f(y))/A$ exactly. The same holds for any noise
density that is constant on $[-2MR,2MR]$.

\emph{Small gaps.} Every transfer function with $\psi(0)=0$ and $\psi'(0)=\lambda>0$
is linear to first order, so the linear link governs the regime $f(x)-f(y)\to0$ that
dominates near the optimum. \Cref{sec:general} makes this quantitative.

\begin{theorem}\label{thm:lin-upper}
Let $Q\subset \B(c, R)^d$, $f\in\Fm$ with $\mu\ge0$ and Assume~\eqref{eq:linear-link}. Then, for any $N\ge 1$, Algorithm \ref{alg:value} guarantees
\[
  \E[f(\hat x_N)]-f^*\ \le\ \min\Big\{MR,\ \frac{C\,d^{3/2}\log^{p}(ed(N+1))}{\lambda\sqrt N}\Big\}
\]
for universal constants $C,p>0$.
\end{theorem}

The proof is given in \cref{proof-thm-lin-upper}. Algorithm~\ref{alg:value} modifies the bandit convex optimization ($\mathsf{BCO}$)
method of \cite{fokkema2024online} and requires neither $\mu$
nor the horizon $N$ as input. The method of \citet{fokkema2024online} is a bandit
version of online Newton step: in every round it samples a point from a Gaussian
$\mathcal N(\mu_t,\Sigma_t)$, rescales it radially into $Q$, and updates $\mu_t$ and
$\Sigma_t$ using an estimate of a quadratic surrogate of the loss. Importantly, Algorithm \ref{alg:value} achieves an optimal rate of $N$, which is stated in the following theorem.
\begin{theorem}\label{thm:lin-lower}
Let $\B^d(c,r)\subseteq Q\subseteq \B^d(c,R)$ and assume ~\eqref{eq:linear-link}.
Then, for every $N\ge1$,
\[
  \cL_N(\Fzero,o)\ \ge\ \min\Big\{\frac{Mr}{8},\ \frac{\sqrt3\,d}{32\lambda\sqrt N}\Big\}.
\]
If $0<\mu\le M/(8R)$, then
\begin{align*}
  \cL_N(\Fm,o)\ \ge\ \max\Big\{&\min\Big\{\frac{21Mr}{256},\ \frac{3\sqrt3\,d}{128\lambda\sqrt N}\Big\},\\
  &\min\Big\{\frac{\mu r^2}{12},\ \frac{d^2}{144\lambda^2\mu R^2N}\Big\}\Big\}.
\end{align*}
\end{theorem}

The proof is given in \cref{proof:thm-lin-lower}.
\paragraph{Discussion.} The two results \cref{thm:lin-upper,thm:lin-lower} show that the optimal exponent is $N^{-1/2}$ for convex
and strongly convex non-smooth objectives, with a $\sqrt d$ gap that mirrors the current
gap between $\widetilde O(d^{3/2}\sqrt N)$ \citep{fokkema2024online} and $\Omega(d\sqrt N)$
\citep{shamir2013complexity} for bandit convex optimization. The hard instances are
 $f_v(x)=\tfrac\mu2\|x\|^2+a\sum_i|x_i-v_is|$, which are strongly convex
yet force the $d/(\lambda\sqrt N)$ term. So, unlike for first-order oracles, strong
convexity does not improve the exponent, in line with the phenomenon identified by
\citet{shamir2013complexity} for noisy zeroth-order feedback. 

Despite that, for each
comparison query, Algorithm \ref{alg:value} solves a convex quadratic subproblem over a subset of $Q$ and maintains and
factorizes a $d\times d$ matrix, which costs $O(d^2)$ to $O(d^3)$ arithmetic operations per
round. It requires placing $Q$ in a special affine position, for which no efficient procedure
is known for general $Q$. In that position, the number of iterations is $d^{3/2}$.
Positions computable in polynomial time give $d^{7/4}$ instead of $d^{3/2}$
\cite{fokkema2024online}, and for symmetric $Q$ the two rates coincide. 
The smooth-case method of \cite{lu2026online} faces the same obstacle: it requires a special convex function, which exists for every convex body but
is not known to be efficiently computable in general. To address this issue, we introduce a rigorous analysis of simpler \Cref{alg:ardg}, which is able to work with general transfer functions.
\begin{figure*}[t]
\centering
\includegraphics[width=0.3\textwidth]{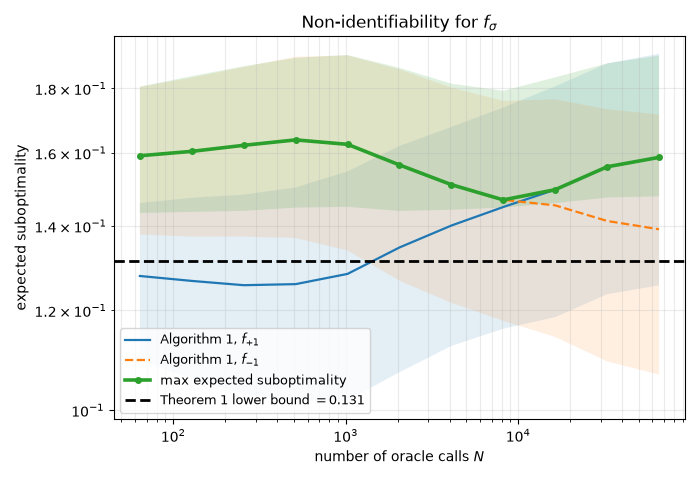}%
\includegraphics[width=0.3\textwidth]{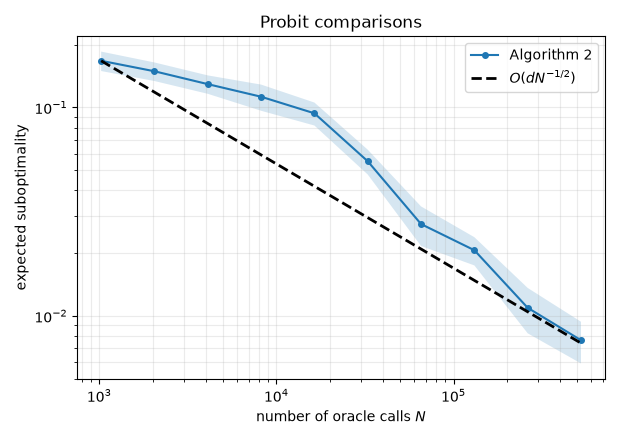}\\[2mm]
\includegraphics[width=0.60\textwidth]{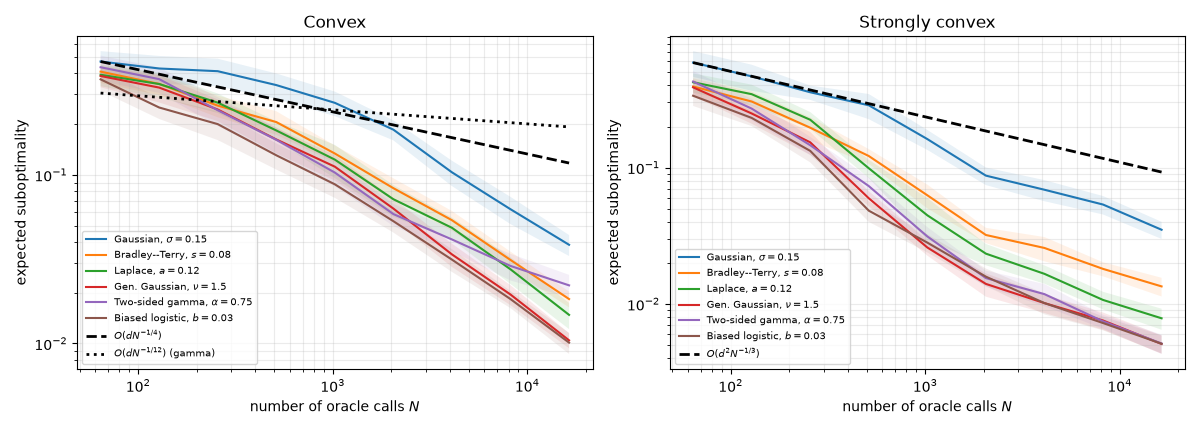}
\caption{Experiments. Upper left: order oracles fail on the instances of \cref{thm:nonident}. Upper right: probit calibration restores convergence. Bottom: ARDG under diverse transfer functions. Shaded regions show $\pm2$ standard errors.}
\label{fig:experiments}
\end{figure*}

\subsection{General Transfer Functions}
\label{sec:general}

For non-linear $\psi$ we analyze \cref{alg:ardg}, an adaptive modification of the
relative gradient descent of \citet{saha2025dueling} that handles non-smooth
objectives, uses $m$ repeated comparisons per iteration, and needs neither the horizon
nor the diameter. At each iteration, it queries several $x_k\pm h_ke_k$ (which can lie out of $Q$) for different $e_k$, averages them and then projects the result onto the domain. As is standard in
two-point methods \citep{flaxman2005online,lu2026online}, we therefore assume that $f$ is
convex and $M$-Lipschitz not only on $Q$, but on $U:=\{x:\operatorname{dist}(x,Q)<2h_{\max}\}$, and
$\mu$-strongly convex there when $\mu>0$.

\begin{algorithm}[t]
\caption{ARDG-$m$: Adaptive Relative Gradient Descent}
\label{alg:ardg}
\begin{algorithmic}[1]
\Require $x_0\in Q$; radii $h_k\in(0,h_{\max}]$; steps $\gamma_k>0$; repetitions $m\ge1$
\For{$k=0,1,2,\dots$}
  \State draw $e_k\sim\Unif(\S^{d-1})$ independently of the past
  \State query $(x_k+h_ke_k,\ x_k-h_ke_k)$ $m$ times; $\bar O_k := $ be the mean response
  \State $x_{k+1}\gets\Projp_Q\big(x_k-\gamma_k\bar O_k e_k\big)$
\EndFor
\State \Return $\bar x_N=\frac1N\sum_{k<N}x_k$ (convex) or
$\bar x_N=\frac{2}{N(N+1)}\sum_{k=1}^N kx_k$ ($\mu>0$)
\end{algorithmic}
\end{algorithm}

Fix any $\lambda>0$ and define the linearization remainder and the range of $\psi$,
\begin{equation}\label{eq:rho}
  \rho(T):=\sup_{|t|\le T}|\psi(t)-\lambda t|,\quad
  \Psi(T):=\sup_{|t|\le T}|\psi(t)|,
\end{equation}
The quantity $\rho$ measures how far
$\psi$ is from its best linear approximation on the gaps that the algorithm actually
probes, which never exceed $2Mh_k$.

\begin{assumption}[Linearizability]\label{as:q}
There are $q,L,T_0>0$ such that $\rho(T)\le LT^{1+q}$ for all $T\le T_0$, and
$2Mh_{\max}\le T_0$.
\end{assumption}

By using these estimates, we can formulate two main theorems, which describe how Algorithm \ref{alg:ardg} is capable of covering various types of noise distributions. 

\begin{theorem}\label{thm:ardg-convex}\footnote{In theorems \ref{thm:ardg-convex} and \ref{thm:ardg-sc}, value of $G_k$ depends only on $m$ and $\Psi$.}
Let $f \in \Fzero(U)$ and $\lambda>0$. 

\textup{(a)} Then, for any $N\ge1$, $(h_k)\subset(0,h_{max}]$ and any $(\gamma_k)$, Algorithm \ref{alg:ardg} produces $\bar x_N$ such that
\begin{align*}
\E f(\bar x_N)-f^*\le&\frac1N\Big[
\frac d{4\lambda}\sum_{k<N}\frac{\gamma_kG_k^2}{h_k}
+M\sum_{k<N}h_k\\&+\frac{D\sqrt d}{2\lambda}\sum_{k<N}\frac{\rho(2Mh_k)}{h_k}+\frac{dD^2}{4\lambda h_{N}\gamma_{N}}\Big].
\end{align*}
\textup{(b)} More precisely, with
\[
  \gamma_k = \frac{\gamma}{\sqrt{k+1}}
  \quad\text{and}\quad
  h_k = \min\{h_{\max},\, c(k+1)^{-\beta}\}
\]
for any \(\gamma, c > 0\):

\begin{itemize}
  \item If \(q \ge 1\) and \(\beta = \frac14\), then
  \[
    \E f(\bar x_N) - f^* = O\!\left(d\,N^{-1/4}\right).
  \]

  \item If \(0 < q < 1\) and \(\beta = \frac{1}{2(1+q)}\), then
  \[
    \E f(\bar x_N) - f^* = O\!\left(d\,N^{-q/(2(1+q))}\right).
  \]
\end{itemize}
\end{theorem}

\begin{theorem}\label{thm:ardg-sc}
Let, in addition, $f$ be $\mu$-strongly convex on $U$, and let
$\gamma_k=(2d)({\lambda\mu h_k(k+1)})^{-1}$.  

\textup{(a)} Then, for any $N\ge1$ and $(h_k)\subset(0,h_{max}]$, Algorithm \ref{alg:ardg} produces $\bar x_N$ such that
\begin{align*}
\E f(\bar x_N)-f^*\le\frac{2}{N(N+1)}\sum_{k=1}^Nk\Big[&
\frac{d^2G_k^2}{2\lambda^2\mu h_k^2(k+1)}+Mh_k\\[-1mm]
&+\frac{d\,\rho(2Mh_k)^2}{4\lambda^2\mu h_k^2}\Big].
\end{align*}
\textup{(b)} More precisely, with
\[
  h_k = \min\{h_{\max},\, c(k+1)^{-\beta}\}
\]

\begin{itemize}
  \item If \(q \ge \frac12\) and \(\beta = \frac13\), then
  \[
    O\!\left(d^2N^{-1/3}\right).
  \]

  \item If \(0 < q < \frac12\) and \(\beta = \frac{1}{2(1+q)}\), then
  \[
    O\!\left(d^2N^{-q/(1+q)}\right).
  \]
\end{itemize}
\end{theorem}

The proof, explicit constants and optimal choices of $c$ and $\gamma$ for both convex and strongly convex regimes are
given in \cref{proof:thm-argd-convex} and \cref{proof:thm-ardg-strongly-convex} respectively. Tuning
$\gamma$ and $c$ in $d$ improves the leading coefficient: in the convex case to
$d^{3/4}$ for $q\ge1$ and $d^{\frac{1+2q}{2(1+q)}}$ for $0<q<1$, and in the strongly convex case
to $d^{4/3}$ for $q\ge\frac12$ and $d^{\frac{1+2q}{1+q}}$ for $0<q<\frac12$. For a linear
$\psi$, these factors become $d^{1/2}$ and $d^{2/3}$, respectively.

All terms in both theorem statements have a clear interpretation: the first is the variance of the estimator, the second the bias of ball smoothing, and the third the departure of the link from linearity. The additional fourth term in Theorem \ref{thm:ardg-convex} is the usual last-iterate term. Choosing the radii is thus a trade-off between bias and variance in which $q$
determines how fast the non-linear bias decays.

\paragraph{Families of allowed noise distributions.}
For $q\ge1$ (convex) and $q\ge\frac12$ (strongly convex), non-linearity is free:
the rates coincide with those for the exactly linear link and only constants change.
In particular, every transfer function that is $C^2$ near $0$ with $\psi(0)=0$ and
$\psi'(0)>0$ has $q\ge1$, and odd $C^3$ functions have $q\ge2$.

For smaller $q$, the rate degrades continuously, with exponents
$\frac q{2(1+q)}$ and $\frac q{1+q}$ that agree with the linear-link exponents at $q=1$
and $q=\frac12$ and vanish as $q\downarrow0$.

For an oracle model $O(x,y)=\text{sign}(f(x)-f(y)+\xi)$,

where $\xi$ has symmetric density $p$, the transfer function is
$\psi(t)=2F_\xi(t)-1$, where $F_\xi$ is the distribution function
of~$\xi$. If $0<p(0)<\infty$ and
$|p(u)-p(0)|\le C|u|^q$ for $|u|\le T_0$, then, for
$0<T\le T_0$,
\begin{equation*}
    \lambda=2p(0),
    \qquad
    \rho(T)\le\frac{2C}{1+q}T^{1+q}.
\end{equation*}

For example, if \(\xi \sim \mathrm{Laplace}(a)\), then \(\lambda = 1/a\) and \(q = 1\).
If each observation carries its own Laplacian noise with possibly different
parameters, i.e.\ \(O(x,y) = \operatorname{sign}\bigl(f(x) + \xi_x - f(y) - \xi_y\bigr)\),
then the resulting noise on the gap \(f(x)-f(y)\) is a convolution and therefore follows a two-sided Gamma distribution. The more general result for generalized Gaussian with arbitrary shape $\nu$ is shown in \Cref{tab:links}, where most of common distributions are listed. All proofs are given in \cref{app:links}.
\begin{table}[t]
\caption{Noise families and their exponents}
\label{tab:links}
\begin{center}\small
\begin{tabular}{@{}lll@{}}
\toprule
Rater noise / model & $\lambda$ & $q$\\
\midrule
Uniform on $[-A,A]$, $A\ge2MR$
    & $A^{-1}$ & $>0$\\
    \addlinespace
Thurstone -- $\cN(0,\sigma^2)$, $\sigma>0$
    & $(\sigma\sqrt\pi)^{-1}$ & $2$\\
        \addlinespace

Bradley--Terry -- Logistic($s, 0$), $s>0$
    & $(2s)^{-1}$ & $2$\\
        \addlinespace

Laplace, $a>0$
    & $a^{-1}$ & $1$\\
        \addlinespace

Generalized Gaussian, $\nu>0$
    & $\frac{\nu}{a\Gamma(1/\nu)}$ & $\nu$\\
        \addlinespace
Two-sided gamma, $\alpha\in(\frac12,1)$
    & see \eqref{eq:gamma-lambda} & $2\alpha-1$\\
        \addlinespace

Logistic($s, b$) $s, b>0$
    & $\frac{1}{2s}\text{sech}^2\!\left(\frac{b}{2s}\right)$ & $1$\\
        \addlinespace

Any $C^2$, $\psi(0)=0<\psi'(0)$
    & $\psi'(0)$ & $\ge1$\\
    
\bottomrule
\end{tabular}
\end{center}
\end{table}
\paragraph{Comparison with prior work.}
The authors of \citet{saha2025dueling} assume that $\psi$ is differentiable, odd, and
sign-preserving, with $\psi'(t)\ge cp|t|^{p-1}$ near $0$, and that $f$ is smooth. In
\citet{lu2026online} it is argued that $\psi$ must be odd for preferences to be
symmetric, and $|\psi''|\le M_2$ is assumed on the attainable range; this is our
condition with $q=1$ and $L=M_2/2$. Our assumption restricts only the local remainder
$\rho$ and covers links outside both settings: two-sided gamma evaluation noise
\citep{jamieson2012query} has $\psi''$ unbounded at the origin; generalized Gaussian
noise with $\nu<1$ gives a $\psi$ that is not $C^2$; restricted raters with
$\psi(t)=\min\{\max\{\lambda t,-a\},a\}$ are not differentiable; and raters with position
bias are not odd (for example, biased logistic, i.e.\ Logistic($s,b$) with $b>0$). The
last case does not even satisfy $\psi(0)=0$, which our method tolerates.
\begin{proposition}[Position bias]\label{prop:bias}
If $\psi(t)=b+\lambda t+\tilde r(t)$ with a constant bias $b\in(-1,1)$, then
\cref{thm:ardg-convex,thm:ardg-sc} hold with $m=1$ and $\rho$ replaced by
$\tilde\rho(T):=\sup_{|t|\le T}|\tilde r(t)|$.
\end{proposition}

\begin{proof}
Since $e_k$ is uniform on the sphere and independent of $\cF_k$,
$\E[b\,e_k\mid\cF_k]=0$, so the bias $b$ leaves the mean of the estimator unchanged;
$\|g_k\|\le1$ is unaffected.
\end{proof}

The decaying radii $(h_k)$ are not an artifact of the analysis. Consider any algorithm whose
$k$-th query is a symmetric pair $x_k\pm h_ke_k$ with any adaptive choice of $x_k$, $e_k$, and $h_k\le H$, subject to the radius budget $\sum_{k<N}h_k^2\le B_N$. Then, the following theorem holds.

\begin{theorem}\label{thm:restricted}
Let $Q=\B(R)^d$, $\psi(t)=\lambda t$, and $H=(2\lambda M)^{-1}$. Then, for $m =1$ and every $N\ge1$,
 $$\mathcal{L}_N(\Fzero(Q), o)\ge\min\{\frac{MR}{16},\frac{dR}{128\lambda\sqrt{B_N}}\}.$$

If $M>2\mu R$ and $\bar M:=\min\{\mu R,\frac{M-2\mu R}2\}$, then $$\mathcal{L}_N(\Fm(\mathbb{B}^d(2R)), o)\ge\min\left\{\frac{\bar M^2}{8\mu},\frac{3d^2}{1024\lambda^2\mu B_N}\right\}.$$
\end{theorem}

The same asymptotic bounds hold for any $m > 0$. These are the schedules of \cref{thm:ardg-convex,thm:ardg-sc}, $B_N=O(N^{1/2})$ and
$B_N=O(N^{1/3})$ respectively, so the lower bounds become $\Omega(N^{-1/4})$ and
$\Omega(N^{-1/3})$. The rates of \cref{alg:ardg} are therefore exact for its radius
schedules, and the only way to improve them is to probe at a non-vanishing radius, which is what presented in Algorithm \ref{alg:value}. The proof and explicit constants in case where $m > 1$, are given in \ref{proof:restricted}.
\section{Experiments}
\label{sec:experiments}

We test the three qualitative predictions of the theory. Details are in
\cref{app:experiments}.
\paragraph{Order oracles fail in general.}
On the instances of \cref{thm:nonident} with $R=1,M=1,d=2,\mu=0.05$, \cref{alg:ardg} with $m=1$ is run
with $30$ random seeds and reaches errors $0.141$ and $0.150$ on $f_\sigma$,
$\sigma\in\{-1,+1\}$, after $2^{16}$ calls — above the predicted floor $0.131$ (\cref{fig:experiments}, upper left).
\paragraph{Probit calibration restores convergence.}
For the normalized probit model with $d=2$, $M=1/2$, $\sigma=1$ and $f(x)=\tfrac{M}{\sqrt d}\|x-x_*\|_1$, the convex configuration of \cref{alg:probit} is run with $10$ random seeds for budgets $N\in\{2^{10},\ldots,2^{19}\}$. The empirical error follows the predicted $O(dN^{-1/2})$ rate (\cref{fig:experiments}, upper right).

\paragraph{ARDG under diverse transfer functions.}
We minimize $f(x)=\tfrac{\mu}{2}\|x-x_*\|_2^2+\tfrac{1}{\sqrt d}\|x-x_*\|_1$ with $d=5$, in the convex ($\mu=0$) and strongly convex ($\mu=1$) cases, for budgets $N\in\{2^{6},\ldots,2^{14}\}$. We consider Gaussian noise with $\sigma=0.15$, Bradley--Terry with $s=0.08$, Laplace
noise with $a=0.12$, generalized Gaussian noise with $\nu=1.5$ and $a=0.12$, two-sided
gamma with $\alpha=0.75$ and $\beta=12$, and biased logistic with $s=0.05$ and $b=0.03$.
Results are averaged over $20$ random seeds. The black dashed lines show the predicted $O(d^{3/4}N^{-1/4})$ and
$O(d^{4/3}N^{-1/3})$ rates; for two-sided gamma in the convex case the
rate is slower (\cref{thm:ardg-convex}), as the thin dashed line shows. All six
methods converge, in several cases faster than the bounds
(\cref{fig:experiments}, bottom). Fitted log-log slopes (second half of the run, in legend order) are
$-0.716,-0.715,-0.750,-0.854,-0.547,-0.777$ (convex) and
$-0.512,-0.526,-0.615,-0.559,-0.627,-0.606$ (strongly convex).
They are at least as steep as the worst-case exponents. The position bias has no visible effect, as \cref{prop:bias}
predicts.

\section{Conclusion}
\label{sec:conclusion}
Comparison feedback suffices for non-smooth convex optimization only when the law of a
comparison is tied to the objective. For Stochastic Order Oracles this tie can fail:
comparisons on a shared realization are invariant to positive rescalings of the sample
loss, so no algorithm guarantees vanishing error, even when every sample loss is convex
and Lipschitz and the mean is strongly convex (\cref{thm:nonident}). Under the probit
calibration~\eqref{eq:probit-model}, the minimax error is
$\min\{MR,\ d\sigma R/\sqrt N,\ d^2\sigma^2/(\mu N)\}$ up to a factor depending only on
$M/\sigma$ (\cref{thm:probit}).

For Preference Oracles with a linear transfer function, the minimax error lies between
$\Omega(d/(\lambda\sqrt N))$ and $\widetilde O(d^{3/2}/(\lambda\sqrt N))$
(\cref{thm:lin-upper,thm:lin-lower}). Unlike in the probit case, strong convexity does
not improve the $N^{-1/2}$ exponent. For general transfer functions, \cref{alg:ardg}
converges at rates determined only by the linearization exponent $q$ of $\psi$ at the
origin, requiring neither oddness, differentiability, nor sign-preservation of $\psi$.
For $q\ge1$ ($q\ge\frac12$ under strong convexity), its rates $N^{-1/4}$ ($N^{-1/3}$)
match the lower bound for methods with the same radius budget
(\cref{thm:ardg-convex,thm:ardg-sc,thm:restricted}).

\paragraph{Limitations and open problems.}
A $\sqrt d$ gap remains for the linear link. Linear-link comparisons and bounded-noise
zeroth-order queries reduce to each other, so closing this gap would also close the
corresponding gap for bandit convex optimization. The method attaining $N^{-1/2}$ is not
computationally efficient for general $Q$. For general transfer functions, the question about optimal method is still open. For Stochastic Order Oracles, it remains open to
characterize the noise models under which comparisons identify the mean objective.
\subsubsection*{AI use statement}
We used AI tools substantially, mainly for formatting and for checking
hypotheses. All mathematical results, proofs, and experiments were verified by the
authors.
\subsubsection*{Acknowledgments}
We want to thank Valentin LePlat for helpful recommendations.
\subsubsection*{Support}
The research was supported by Russian Science Foundation (project No. 21-71-30005-\text{П}), https://rscf.ru/en/project/21-71-30005/.

\clearpage

\appendix
\thispagestyle{empty}
\onecolumn
\begin{center}{\Large\bfseries Supplementary Material}\end{center}

\newcommand{\hmax}{h_{max}}
\newcommand{\norm}[1]{\|#1\|}

\newcommand{\Fk}{\mathcal{F}_k}

\section{Divergence of Stochastic Order Oracle}
\subsection{Proof of \Cref{thm:nonident}}\label{app:nonident}

Let $\xi\in \{\pm1\}$ be a random variable, such that $\Pr(\xi=1)=\Pr(\xi=-1)=\tfrac12$.
For $\sigma\in\{+1,-1\}$ and $x\in Q=\B^d(R)$ define
\[
  f_\sigma(x,\xi)=\big(1+\tfrac34\sigma\xi\big)\,q(x,\xi),\qquad
  q(x,\xi)=\frac\mu2\|x\|^2+\frac M4\,\xi x_1+\frac M{32}|x_1|.
\]
Note that $1+\frac34\sigma\xi\in\{\tfrac14,\tfrac74\}$ for all $\sigma,\xi$.

We will separate our proof to steps, so it will be simpler to understand.

\smallskip
\noindent\emph{Step 1.}
Fix $\xi\in\{-1,+1\}$. The function $q(\cdot,\xi)$ is the sum of a convex
quadratic, a linear function and the convex function $x\mapsto|x_1|$, so it is
convex; since the factor $1+\tfrac34\sigma\xi$ is positive, $f_\sigma(\cdot,\xi)$
is convex as well. 

For $x,y\in Q$ we have
$\big|\tfrac\mu2\|x\|^2-\tfrac\mu2\|y\|^2\big|\le\mu R\|x-y\|$, and the maps
$x\mapsto x_1$ and $x\mapsto|x_1|$ are $1$-Lipschitz. Hence $q(\cdot,\xi)$ is
Lipschitz on $Q$ with constant
\[
  \mu R+\frac M4+\frac M{32}=\mu R+\frac{9M}{32}\le\frac{13M}{32},
\]
where we used $\mu R\le M/8$. Multiplying by $1+\tfrac34\sigma\xi\le\tfrac74$,
the sample loss $f_\sigma(\cdot,\xi)$ is Lipschitz with constant at most
$\tfrac74\cdot\tfrac{13M}{32}=\tfrac{91M}{128}\le M$. Thus,
$f_\sigma(\cdot,\xi)\in\Fzero(Q)$ for every $\xi$. (They are even strongly convex, but it is excessive).

\smallskip
\noindent\emph{Step 2.}
Since $\E\big[1+\tfrac34\sigma\xi\big]=1$ and
$\E\big[\xi\big(1+\tfrac34\sigma\xi\big)\big]=\E\xi+\tfrac34\sigma\E\xi^2
=\tfrac34\sigma$,
\[
  f_\sigma(x)=\E_\xi\big[f_\sigma(x,\xi)\big]
  =\frac\mu2\|x\|^2+\frac{3\sigma M}{16}\,x_1+\frac M{32}|x_1|.
\]
This function is $\mu$-strongly convex and, by the same argument as in Step~1,
Lipschitz on $Q$ with constant
$\mu R+\tfrac{3M}{16}+\tfrac M{32}=\mu R+\tfrac{7M}{32}\le\tfrac{11M}{32}\le M$.
Hence $f_\sigma\in\Fm(Q)$.

Let $x\in Q$ and $r=\|x\|\le R$. Using $\sigma x_1\ge-|x_1|$ and $|x_1|\le r$,
\[
  f_\sigma(x)\;\ge\;\tfrac\mu2 r^2-\tfrac{3M}{16}|x_1|+\tfrac M{32}|x_1|
  \;=\;\tfrac\mu2 r^2-\tfrac{5M}{32}|x_1|
  \;\ge\;\tfrac\mu2 r^2-\tfrac{5M}{32}\,r=:g(r).
\]
For $r\in[0,R]$, $g'(r)=\mu r-\tfrac{5M}{32}\le\mu R-\tfrac{5M}{32}
\le-\tfrac M{32}<0$, so $g$ attains its minimum over $[0,R]$ at $r=R$. At
$x=-\sigma Re_1$ we have $\|x\|=R$, $\sigma x_1=-|x_1|$ and $|x_1|=R$, so all
inequalities above hold with equality. Therefore
\begin{equation}\label{eq:nonident-min}
  x_\sigma^*=-\sigma Re_1,\qquad
  f_\sigma^*=\frac{\mu R^2}{2}-\frac{5MR}{32}.
\end{equation}
The minimizers $x_+^*=-Re_1$ and $x_-^*=Re_1$ lie at distance $2R$.

\smallskip
\noindent\emph{Step 3.}
For all $x,y\in Q$ and $\xi\in\{-1,+1\}$,
\[
  f_\sigma(x,\xi)-f_\sigma(y,\xi)
  =\big(1+\tfrac34\sigma\xi\big)\big(q(x,\xi)-q(y,\xi)\big),
\]
and the factor $1+\tfrac34\sigma\xi$ is positive. Hence
\begin{equation}\label{eq:sign-invariance}
  \sign\big(f_+(x,\xi)-f_+(y,\xi)\big)
  =\sign\big(f_-(x,\xi)-f_-(y,\xi)\big)
  =\sign\big(q(x,\xi)-q(y,\xi)\big),
\end{equation}
and in particular the two differences vanish simultaneously.

Fix $N\ge0$ and an algorithm. If the algorithm makes fewer than $N$ calls, we
pad it with calls whose answers it ignores, so we may assume that it makes
exactly $N$ calls. Let $U$ be its internal randomness, $\xi_1,\dots,\xi_N$ the
oracle samples, and $b_1,\dots,b_N$ the tie-breaking coins of
\cref{def:order}. These variables are mutually independent, and their joint law
does not depend on $\sigma$. Write $\omega=(U,\xi_{1:N},b_{1:N})$. The algorithm
is described by measurable maps $\Phi_1,\dots,\Phi_N,\Psi$, which do not depend
on $\sigma$, such that
\[
  (X_t,Y_t)=\Phi_t\big(U,O_{1:t-1}\big)\quad(t\le N),\qquad
  \hat x_N=\Psi\big(U,O_{1:N}\big).
\]
Under instance $\sigma$ the oracle returns
$O_t^\sigma=\phi\big(f_\sigma(X_t,\xi_t)-f_\sigma(Y_t,\xi_t),\,b_t\big)$, where
$\phi(s,b)=\sign(s)$ if $s\ne0$ and $\phi(0,b)=b$.

We claim that $O_t^+(\omega)=O_t^-(\omega)$ for all $\omega$ and all $t\le N$.
For $t=1$, the query $(X_1,Y_1)=\Phi_1(U)$ is the same for both instances, and
\eqref{eq:sign-invariance} with $\xi=\xi_1$ gives $O_1^+=O_1^-$. If the claim
holds for all indices smaller than $t$, then
$(X_t,Y_t)=\Phi_t(U,O_{1:t-1})$ is again the same for both instances, and
\eqref{eq:sign-invariance} with $\xi=\xi_t$ gives $O_t^+=O_t^-$. By induction,
the two transcripts coincide, so the output $\hat x_N=\Psi(U,O_{1:N})$ is one and
the same function $\hat x_N(\omega)$ for $\sigma=+1$ and $\sigma=-1$. Because the
law of $\omega$ does not depend on $\sigma$, the output has the same
distribution under both instances.

\smallskip
\noindent\emph{Step 4.}
Let $L_\sigma=\E_\sigma\big[f_\sigma(\hat x_N)-f_\sigma^*\big]$. By Step~3,
$L_\sigma=\E\big[f_\sigma(\hat x_N(\omega))\big]-f_\sigma^*$, with the same
random point $\hat x_N(\omega)$ for both values of $\sigma$. The linear terms of
$f_+$ and $f_-$ cancel in their average:
\[
  \tfrac12\big(f_+(x)+f_-(x)\big)=\tfrac\mu2\|x\|^2+\tfrac M{32}|x_1|\ \ge\ 0
  \qquad\text{for all }x\in Q.
\]
Combining this with \eqref{eq:nonident-min},
\begin{align*}
  \max\{L_+,L_-\}
  \;\ge\;\tfrac12\,(L_++L_-)
  &=\E\Big[\tfrac12\big(f_+(\hat x_N)+f_-(\hat x_N)\big)\Big]
    -\tfrac12\big(f_+^*+f_-^*\big)\\
  &\ge-\tfrac12\big(f_+^*+f_-^*\big)
   \;=\;\frac{5MR}{32}-\frac{\mu R^2}{2}.
\end{align*}
Finally, $\mu\le M/(8R)$ gives $\mu R^2/2\le MR/16$, hence
$\max\{L_+,L_-\}\ge\tfrac{5MR}{32}-\tfrac{2MR}{32}=\tfrac{3MR}{32}$. The
bound holds for every $N$ and the instances depend only on $(d,R,M,\mu)$, so no
algorithm using this oracle can make the error vanish on both instances.

\section{Probit Stochastic Comparison Oracle}

\subsection{Proof of \Cref{thm:probit}(a)}\label{proof:thm_probit}

\emph{Step 1.}
Let $q(t):=\Phi^{-1}\big(\tfrac{1+t}2\big)$ for $t\in(-1,1)$. For a pair $y\ne z$ set
$\nu:=\E[O\mid y,z]=2\Phi\big(\frac{f(y)-f(z)}{\sigma\|y-z\|}\big)-1$. Then
\begin{equation}\label{eq:inverse-link}
  q(\nu)=\frac{f(y)-f(z)}{\sigma\|y-z\|},
  \qquad
  |\nu|\le m_0:=2\Phi(M/\sigma)-1<1,
\end{equation}
using $|f(y)-f(z)|\le M\|y-z\|$. By \cref{lem:probit-taylor}, $q$ has the power series
\eqref{eq:probit-series} with nonnegative coefficients satisfying \eqref{eq:probit-bound}
and the recursion \eqref{eq:probit-coefficients2}.  Below, $r_1,r_2,\theta,k,p_j,V,u,B_0$ are the constants
of \cref{lem:probit-estimator} for this $m_0$.

\emph{Step 2.}
Fix $h\in(0,r/2]$ and let $K_h:=(1-\frac hr)Q$. If $x=(1-\frac hr)y$ with $y\in Q$ and
$\|e\|\le1$, then
\[
  x+he=\big(1-\tfrac hr\big)y+\tfrac hr\,(re)
\]
is a convex combination of $y\in Q$ and $re\in \B^d(r)\subseteq Q$, so $x+he\in Q$. Thus
all queries $x\pm he$ with $x\in K_h$, $e\in\S^{d-1}$ are feasible, and
$f_h(x):=\E_{U\sim\Unif(\B^d)}f(x+hU)$ evaluates $f$ only inside $Q$. By
\cref{lem:smoothing} (applied to the Lipschitz extension
$\tilde f(x)=\inf_{y\in Q}\{f(y)+M\|x-y\|\}$, which coincides with $f$ on $Q$), on $K_h$
\begin{equation}\label{eq:probit-smoothing}
  f\le f_h\le f+Mh,\qquad f_h\ \text{is $\mu$-strongly convex},\qquad
  \nabla f_h(x)=\frac d{2h}\,\E_e\big[(f(x+he)-f(x-he))\,e\big],
\end{equation}
with $e\sim\Unif(\S^{d-1})$.

For $x\in K_h$, draw $e\sim\Unif(\S^{d-1})$, compute $T$ by \cref{lem:probit-estimator}
at the pair $(x+he,x-he)$, and set
\begin{equation}\label{eq:poiu}
    \hat g:=d\,\sigma\,T\,e.
\end{equation}
The two points are at distance $2h$, so \eqref{eq:inverse-link} gives
$\E[T\mid e]=\frac{f(x+he)-f(x-he)}{2\sigma h}$, and therefore
\begin{equation}\label{eq:probit-gradient-moments2}
  \E\hat g=\frac d{2h}\,\E_e\big[(f(x+he)-f(x-he))\,e\big]=\nabla f_h(x),
  \qquad
  \E\|\hat g\|^2=d^2\sigma^2\,\E T^2\le d^2\sigma^2V =: G^2.
\end{equation}

\emph{Step 3.}
Given $n\ge1$, set $x_1:=0$ and $x_{t+1}:=\Projp_{K_h}(x_t-\eta_t\hat g_t)$ for
$t=1,\dots,n$, where $\hat g_t$ is defined as in \eqref{eq:poiu} at $x_t$ with fresh
randomness, and output $z:=\sum_{t\le n}w_tx_t/\sum_{t\le n}w_t$. We use two
configurations:
\begin{equation}\label{eq:probit-configurations2}
\begin{aligned}
\mathrm c:&\quad h=\min\{\tfrac r2,\tfrac{RG}{\Gamma\sqrt n}\},\quad \eta_t=\tfrac{R}{G\sqrt n},\quad w_t=1;\\
\mathrm s:&\quad h=\min\{\tfrac r2,\tfrac{G^2}{\mu\Gamma(n+1)}\},\quad \eta_t=\tfrac{2}{\mu(t+1)},\quad w_t=t\qquad(\mu>0),
\end{aligned}
\end{equation}
with $\Gamma:=M(1+R/r)$. Conditionally on the history before step $t$,
\eqref{eq:probit-gradient-moments2} holds at $x_t$.

Let $x^*\in\arg\min_{x\in Q}f$ and $z_h:=(1-\frac hr)x^*\in K_h$. Since
$\|z_h-x^*\|=\frac hr\|x^*\|\le\frac{hR}r$,
\[
  f_h(z_h)-f^*\le f(z_h)+Mh-f(x^*)\le Mh+\frac{MhR}r=\Gamma h .
\]
Let $D_t:=\E\|x_t-z_h\|^2$. As $\Projp_{K_h}$ is nonexpansive and fixes $z_h$,
\[
  \|x_{t+1}-z_h\|^2\le\|x_t-z_h\|^2-2\eta_t\ip{\hat g_t}{x_t-z_h}+\eta_t^2\|\hat g_t\|^2 .
\]
By strong convexity of $f_h$ and \eqref{eq:probit-smoothing},
\[
  \ip{\nabla f_h(x_t)}{x_t-z_h}\ge f_h(x_t)-f_h(z_h)+\tfrac\mu2\|x_t-z_h\|^2
  \ge f(x_t)-f^*-\Gamma h+\tfrac\mu2\|x_t-z_h\|^2 .
\]
Taking expectations and using \eqref{eq:probit-gradient-moments2},
\begin{equation}\label{eq:probit-recursion2}
  \E f(x_t)-f^*\le\frac{D_t-D_{t+1}}{2\eta_t}-\frac\mu2D_t+\frac{\eta_tG^2}2+\Gamma h .
\end{equation}
\emph{Configuration $\mathrm c$.} Drop $-\frac\mu2D_t\equiv0$, sum
\eqref{eq:probit-recursion2} over $t\le n$, use $D_1=\|z_h\|^2\le R^2$, and apply
Jensen's inequality:
\[
  \E f(z)-f^*\le\frac{R^2}{2\eta n}+\frac{\eta G^2}2+\Gamma h
  =\frac{RG}{\sqrt n}+\Gamma h\le\frac{2RG}{\sqrt n}.
\]
\emph{Configuration $\mathrm s$.} Multiply \eqref{eq:probit-recursion2} by $t$. With
$\eta_t=\frac2{\mu(t+1)}$ the distance terms become
\[
  \frac{\mu t(t+1)}4(D_t-D_{t+1})-\frac{\mu t}2D_t
  =\frac\mu4\big(t(t-1)D_t-(t+1)t\,D_{t+1}\big),
\]
which telescope to $-\frac\mu4n(n+1)D_{n+1}\le0$. Moreover
$\sum_{t\le n}t\,\frac{\eta_tG^2}2=\frac{G^2}\mu\sum_{t\le n}\frac t{t+1}\le\frac{nG^2}\mu$.
Dividing by $\sum_{t\le n}t=\frac{n(n+1)}2$ and applying Jensen's inequality,
\[
  \E f(z)-f^*\le\frac{2G^2}{\mu(n+1)}+\Gamma h\le\frac{3G^2}{\mu(n+1)}.
\]

\emph{Step 4.}
Step~3 analyzes SGD in which update $t$ uses as many calls as the estimator of
\cref{lem:probit-estimator} needs, namely $K_t:=kJ_t$, which is random and unbounded.
\Cref{alg:probit} has a fixed budget of $b$ calls per block, so we must show that
respecting the budget costs little. Fix $b\ge1$ and put
\[
  C_0:=\Big\lceil\max\Big\{2,\frac{2\log B_0}u\Big\}\Big\rceil,\qquad n:=\lfloor b/C_0\rfloor .
\]
If $n=0$, the block outputs $0$, whose error is at most $MR$, because
$f(0)-f^*\le M\|0-x^*\|\le MR$. Assume from now on $n\ge1$, and let $h,\eta_t,w_t$ be
the parameters of the chosen configuration \emph{c} or \emph{s} from
\eqref{eq:probit-configurations2}. They depend on $n$ but not on the random calls.

We compare two procedures driven by the same randomness. For each $t\le n$ let
$e_t\sim\Unif(\S^{d-1})$, $J_t$ with $\Pr(J_t=j)=p_j$, and an infinite sequence of
fresh oracle answers at the pair $(x_t+he_t,\,x_t-he_t)$ be independent of everything
else. Update $t$ reads the first $K_t=kJ_t$ answers of its sequence.
\begin{itemize}
\item The first type of iterations (the \emph{uncapped variant}, since it uses as many
  oracle calls as needed) is the process of Step~3: it performs all $n$ updates, no
  matter how many calls they take in total, and returns
  $z_n=\sum_{t\le n}w_tx_t/\sum_{t\le n}w_t$. It is not a real algorithm, since it may
  use more than $b$ calls; we use it only for the analysis. Step~3 gives
  \begin{equation}\label{eq:uncapped-bound}
    \E f(z_n)-f^*\le\begin{cases}\dfrac{2RG}{\sqrt n}&\text{for }\mathrm c,\\[2mm]
    \dfrac{3G^2}{\mu(n+1)}&\text{for }\mathrm s.\end{cases}
  \end{equation}
\item The \emph{capped variant} (constrained from above on the number of oracle calls)
  is what \cref{alg:probit} does. Before update $t$, it draws $J_t$, and if the calls
  already used plus $kJ_t$ exceed $b$, it stops and returns $0$. Otherwise it performs
  update $t$ exactly as above. After the last update, it returns $z_n$, and in all
  cases it makes dummy calls up to a total of exactly $b$.
\end{itemize}

They differ only on one event. Let $E:=\{\sum_{t\le n}K_t\le b\}$. Since all $K_t\ge1$,
the capped block stops at some update if and only if the partial sums
$\sum_{s\le t}K_s$ exceed $b$ for some $t\le n$, i.e.\ if and only if $E^c$ occurs. So
on $E$ the capped block never stops and returns the same point $z_n$ as the uncapped
process, while on $E^c$ it returns $0$, with error at most $MR$. Since
$z_n\in K_h\subseteq Q$, we have $f(z)-f^*\ge0$ on every outcome. Therefore, writing
$z^{\rm cap}$ for the output of the capped block,
\[
  \E f(z^{\rm cap}_n)-f^*
  =\E\big[(f(z_n)-f^*)\mathbf 1_E\big]+\E\big[(f(0)-f^*)\mathbf 1_{E^c}\big]
  \le\E f(z_n)-f^*+MR\,\Pr(E^c).
\]

Step~3 cannot be applied to the capped block directly: it stops depending on the values
of $J_t$, so its gradient estimates are no longer unbiased. We therefore apply Step~3
to the uncapped process, where unbiasedness holds exactly, and pay for the difference
only on the rare event $E^c$.

The $K_t$ are i.i.d.\ and $\E e^{uK_t}=B_0$ by \cref{lem:probit-estimator}(iii). By
Markov's inequality applied to $e^{u\sum_{t\le n}K_t}$,
\[
  \Pr(E^c)=\Pr\Big(\sum_{t\le n}K_t>b\Big)\le e^{-ub}\,\E e^{u\sum_tK_t}=e^{-ub}B_0^{\,n}\le e^{-ub/2}.
\]
The last step uses $n\log B_0\le\frac b{C_0}\log B_0\le\frac{ub}2$, which is exactly why
$C_0\ge\frac{2\log B_0}u$ was chosen.

Combining \eqref{eq:uncapped-bound}, the bound on $\E f(z^{\rm cap}_n)-f^*$ and the
estimate $\Pr(E^c)\le e^{-ub/2}$, the capped block with configurations \emph{c} and
\emph{s} has error at most
\begin{equation}\label{eq:probit-selection-capped}
U_{\mathrm c}(b)=\frac{2RG}{\sqrt n}+MRe^{-ub/2},\qquad
U_{\mathrm s}(b)=\frac{3G^2}{\mu(n+1)}+MRe^{-ub/2},
\end{equation}
respectively, while the trivial option ``always output $0$'' has error at most
$U_0(b)=MR$. Here $U_{\mathrm s}$ is defined only if $\mu>0$. Therefore,
\cref{alg:probit} computes them before the block starts and runs the option with the
smallest one. Hence, the output $z_b$ of a block with budget $b$ satisfies
\begin{equation}\label{eq:block-guarantee}
  \E f(z_b)-f^*\le\min\{U_0(b),U_{\mathrm c}(b),U_{\mathrm s}(b)\}.
\end{equation}

\emph{Step 5.}
Block $\ell=0,1,2,\dots$ has budget $2^\ell$, restarts from $x_1=0$ with fresh
randomness, and publishes its output when it ends. Fix $N\ge1$ and let
$L:=\lfloor\log_2(N+1)\rfloor-1\ge0$. Blocks $0,\dots,L$ use $2^{L+1}-1\le N$ calls,
while blocks $0,\dots,L+1$ would need $2^{L+2}-1>N$. So, after $N$ calls the published
point $\hat x_N$ is the output of the block with budget
\[
  b=2^L\in\big(\tfrac{N+1}4,\tfrac{N+1}2\big],
\]
and \eqref{eq:block-guarantee} applies to it. We compare its right-hand side with
$\mathcal{T}_\mu(N,R,\sigma)$ using two elementary facts:
$e^{-ub/2}\le\min\{(ub)^{-1/2},(ub)^{-1}\}$ (from $ye^{-y}\le e^{-1}$), and, if $\mu>0$,
$MR\le\frac{M^2}{8\mu}$ (from $\mu R\le M/8$).

\emph{Case $b\ge C_0$.} Then $n\ge1$, $n\ge\frac b{2C_0}$ and $n+1\ge\frac b{C_0}$. With
$G=d\sigma K_0$ and $d\ge1$,
\begin{align*}
  U_{\mathrm c}(b)&\le\Big(2\sqrt2\,K_0\sqrt{C_0}+\frac{M/\sigma}{\sqrt u}\Big)\frac{d\sigma R}{\sqrt b},\\
  U_{\mathrm s}(b)&\le\Big(3C_0K_0^2+\frac{(M/\sigma)^2}{8u}\Big)\frac{d^2\sigma^2}{\mu b},
\end{align*}
and $b>N/4$ gives $b^{-1/2}<2N^{-1/2}$ and $b^{-1}<4N^{-1}$.

\emph{Case $b<C_0$.} Then $n=0$, $\hat x_N=0$, and the error is at most $MR$. Since
$N<4b<4C_0$,
\[
  MR\le2\sqrt{C_0}\,\frac M\sigma\cdot\frac{d\sigma R}{\sqrt N},\qquad
  MR\le\frac{M^2}{8\mu}\le\frac{C_0}2\Big(\frac M\sigma\Big)^2\frac{d^2\sigma^2}{\mu N}.
\]
In both cases the error is at most $C(M/\sigma)$ times each of the three terms of
$\mathcal{T}_\mu(N,R,\sigma)$, where
\[
  C(M/\sigma):=\max\Big\{4\sqrt2\,K_0\sqrt{C_0}+\frac{2M/\sigma}{\sqrt u},\
  12C_0K_0^2+\frac{(M/\sigma)^2}{2u},\ 2\sqrt{C_0}\,\frac M\sigma,\ \frac{C_0}2\Big(\frac M\sigma\Big)^2\Big\}\ \ (\ge1).
\]
Hence $\E f(\hat x_N)-f^*\le C(M/\sigma)\,\mathcal{T}_\mu(N,R, \sigma)$. A simpler
admissible choice is $C(M/\sigma)=64C_0(1+K_0+M/\sigma)^2(1+u^{-1})$, which completes
the proof of part (a).

\subsection{Proof of \Cref{thm:probit}(b)}\label{proof:thm-probit-b}

\emph{Step 1}. Let $\alpha:=\min\{\frac M{4\sqrt d},\frac{\sigma\sqrt d}{4\sqrt N}\}$ and
$a:=\alpha\sqrt d\le\frac M4$. For $v\in\{-1,1\}^d$ let
\[
  f_v(x):=\frac\mu2\|x\|^2-\alpha\ip vx+\frac M2(\|x\|-r)_+ ,
\]
with answers $\sign(f_v(y,Z)-f_v(z,Z))$, $f_v(x,Z)=f_v(x)+\ip Zx$, $Z\sim\cN(0,\sigma^2I_d)$;
these realize \eqref{eq:probit-model} for $f=f_v$. 

Each $f_v$ is $\mu$-strongly convex, since $\frac\mu2\|x\|^2$ is $\mu$-strongly convex, while $-\alpha\ip vx$ is linear and $\frac M2(\|x\|-r)_+$ is convex. Moreover, on $Q$ the three terms are Lipschitz with constants $\mu R$, $\alpha\|v\|=a$ and $\frac M2$ respectively, so the Lipschitz constant of $f_v$ on $Q$ is at most $\mu R+a+\frac M2\le\frac{7M}8$, where $\mu R\le\frac M8$ and $a\le\frac M4$. Hence $f_v\in\Fm(Q)$.

\emph{Step 2}. For a unit vector $u$ and $t>r$,
\[
  f_v(tu)-f_v(ru)=\frac\mu2(t^2-r^2)-\alpha(t-r)\ip vu+\frac M2(t-r)
  \ge(t-r)\Big(\frac M2-a\Big)\ge0 .
\]
Hence $f_v(x)\ge f_v(y)$ for $x\in Q$ and $y:=\Projp_{\B^d(r)}(x)$, and
$\min_{x\in Q}f_v=\min_{x\in\B^d(r)}f_v$.

\emph{Step 3}. On $\B^d(r)$, $f_v(y)=\frac\mu2\|y\|^2-\alpha\ip vy$. Let
$r_*:=r$ if $\mu=0$ and $r_*:=\min\{r,a/\mu\}$ if $\mu>0$, and let
$x_v^*:=\frac{r_*}{\sqrt d}v$, so that $\|x_v^*\|=r_*$ and $\alpha v=\frac a{r_*}x_v^*$.
Expanding both sides gives, for $\|y\|\le r$,
\begin{equation}\label{eq:probit-growth}
  f_v(y)-f_v(x_v^*)-\frac a{2r_*}\|y-x_v^*\|^2
  =\frac12\Big(\frac a{r_*}-\mu\Big)\big(r_*^2-\|y\|^2\big)\ \ge0 .
\end{equation}
The right-hand side is nonnegative: $\frac a{r_*}\ge\mu$ by the choice of $r_*$; if
$r_*=r$ then $\|y\|\le r_*$, and if $r_*=a/\mu<r$ the first factor vanishes. Thus,
$x_v^*$ minimizes $f_v$. Let $\hat v_i:=\sign(y_i)$ (with $\sign0:=1$). A wrong sign
$\hat v_i\ne v_i$ forces $|y_i-x^*_{v,i}|\ge|x^*_{v,i}|=r_*/\sqrt d$, so for every
$x\in Q$
\begin{equation}\label{eq:probit-growth-Q}
  f_v(x)-f_v^*\ \ge\ \frac a{2r_*}\|y-x_v^*\|^2\ \ge\ \frac a{2r_*}\cdot\frac{r_*^2}d\cdot|\{i\in \{1,...,d\}:\hat v_i\ne v_i\}|.
\end{equation}

\emph{Step 4}.
Draw $v$ uniformly from $\{-1,1\}^d$ and fix any algorithm that makes at most $N$ calls
and outputs $\hat x\in Q$. Its output depends on the answers and on its own randomness
$\omega$, which is independent of $v$. For each $v\in\{-1,1\}^d$ let $\mathbb P_v$ be the
joint law of $(\omega,O_1,\dots,O_N)$ when the answers come from $f_v$. Let
$v^{(i)}$ be $v$ with the $i$-th sign flipped.

Fix a queried pair $(y,z)$, $y\ne z$. The answer is the sign of a $\cN(s_v,1)$ variable
with $s_v:=\frac{f_v(y)-f_v(z)}{\sigma\|y-z\|}$. Since the sign is a function of this
Gaussian, the data processing inequality gives
\[
  \KL(\mathbb P_v\|\mathbb P_{v^{(i)}})
  \le\KL\big(\cN(s_v,1)\,\big\|\,\cN(s_{v^{(i)}},1)\big)
  =\tfrac12(s_v-s_{v^{(i)}})^2 .
\]
Now $f_v$ and $f_{v^{(i)}}$ differ only in the linear term, so
\[
  \big(f_v(y)-f_v(z)\big)-\big(f_{v^{(i)}}(y)-f_{v^{(i)}}(z)\big)
  =-\alpha\ip{v-v^{(i)}}{y-z}=-2\alpha v_i(y_i-z_i),
\]
and therefore
\begin{equation}\label{eq:porbit-kl}
  \sum_{i=1}^d\tfrac12(s_v-s_{v^{(i)}})^2
  =\sum_{i=1}^d\frac{2\alpha^2(y_i-z_i)^2}{\sigma^2\|y-z\|^2}
  =\frac{2\alpha^2}{\sigma^2}.    
\end{equation}

\emph{Step 5}.
Let $\hat y:=\Projp_{\B^d(r)}(\hat x)$, $\hat v_i:=\sign(\hat y_i)$, and
$A_{v,i}:=\{\hat v_i\ne v_i\}$. The event $A_{v,i}$ is determined by
$(\omega,O_1,\dots,O_N)$, and $A_{v,i}^c=A_{v^{(i)},i}$ because
$v^{(i)}_i=-v_i$. Hence
\[
  \mathbb P_v(A_{v,i})+\mathbb P_{v^{(i)}}(A_{v^{(i)},i})
  =\mathbb P_v(A_{v,i})+\mathbb P_{v^{(i)}}(A_{v,i}^c)
  \ge1-\operatorname{TV}\big(\mathbb P_v,\mathbb P_{v^{(i)}}\big).
\]
Since $v\mapsto v^{(i)}$ is an involution of $\{-1,1\}^d$, summing over $v$ gives
\[
  \sum_v\mathbb P_v(A_{v,i})
  \ge\frac12\sum_v\Big(1-\operatorname{TV}\big(\mathbb P_v,\mathbb P_{v^{(i)}}\big)\Big).
\]
Let $W:=|\{i:\hat v_i\ne v_i\}|$. Averaging over $v$ and summing over $i$,
\[
  \E_v\E_{\mathbb P_v}W
  \ge\frac12\Big(d-\sum_{i=1}^d\E_v\operatorname{TV}\big(\mathbb P_v,\mathbb P_{v^{(i)}}\big)\Big).
\]
By Pinsker's inequality, the Cauchy--Schwarz inequality over the pairs $(i,v)$, and
\eqref{eq:porbit-kl},
\[
  \sum_{i=1}^d\E_v\operatorname{TV}\big(\mathbb P_v,\mathbb P_{v^{(i)}}\big)
  \le\sum_{i=1}^d\E_v\sqrt{\tfrac12\operatorname{KL}\big(\mathbb P_v\|\mathbb P_{v^{(i)}}\big)}
  \le\sqrt{\frac d2\,\E_v\sum_{i=1}^d\operatorname{KL}\big(\mathbb P_v\|\mathbb P_{v^{(i)}}\big)}
  \le\frac{\alpha\sqrt{dN}}\sigma\le\frac d4,
\]
where the last step uses $\alpha\le\frac{\sigma\sqrt d}{4\sqrt N}$. Thus, $\E_v\E_{\mathbb P_v}W\ge\frac{3d}8$.

\emph{Step 6}. By \eqref{eq:probit-growth-Q}, $f_v(\hat x)-f_v^*\ge\frac{a r_*}{2d}\,W$, so
\[
  \E_v\Big[\E_{\mathbb P_v}f_v(\hat x)-f_v^*\Big]\ge\frac{a r_*}{2d}\cdot\frac{3d}8=\frac{3\,a r_*}{16},
\]
and hence for some $v$ the expected suboptimality of the algorithm on $f_v\in\Fm(Q)$ is at least
$\frac{3ar_*}{16}$. It remains to bound $ar_*$ from below. Recall
$a=\min\{\frac M4,\frac{\sigma d}{4\sqrt N}\}$. If $\mu=0$, then $r_*=r$ and
$ar_*=\frac14\min\{Mr,\frac{d\sigma r}{\sqrt N}\}\ge\frac1{16}\mathcal T_\mu(N,r,\sigma)$.
If $\mu>0$, then $ar_*=\min\{ar,\,a^2/\mu\}$, where $ar\ge\frac14\min\{Mr,\frac{d\sigma r}{\sqrt N}\}$ and
\[
  \frac{a^2}\mu=\frac1{16\mu}\min\Big\{M^2,\frac{d^2\sigma^2}N\Big\}\ge\frac1{16}\min\Big\{8Mr,\frac{d^2\sigma^2}{\mu N}\Big\},
\]
because $\mu\le\frac M{8R}$ gives $\frac{M^2}\mu\ge8MR\ge8Mr$. In both cases $ar_*\ge\frac1{16}\mathcal T_\mu(N,r,\sigma)$, so
\[
  \cL_N(\Fm,O_t)\ge\frac3{16}\cdot\frac1{16}\,\mathcal T_\mu(N,r,\sigma)=\frac3{256}\,\mathcal T_\mu(N,r,\sigma)\ge\frac1{100}\,\mathcal T_\mu(N,r,\sigma),
\]
which completes the proof for block (b).

\subsection{Algorithm pseudo-code}
\paragraph{Constants of \cref{alg:probit}.} For reference:
\begin{equation}\label{eq:probit-constants}
\begin{gathered}
m_0=2\Phi(M/\sigma)-1,\quad r_1=\tfrac{1+m_0}2,\quad r_2=\tfrac{1+r_1}2,\quad \theta=\tfrac{r_1}{r_2},\quad
k=\Big\lceil\tfrac{1-m_0^2}{r_1^2-m_0^2}\Big\rceil,\quad p_j=(1-\theta)\theta^{j-1},\\
V=\tfrac{q(r_2)^2\theta^2}{(1-\theta)^2},\quad u=\tfrac{\log(1/\theta)}{2k},\quad B_0=1+\theta^{-1/2},\quad
C_0=\Big\lceil\max\Big\{2,\tfrac{2\log B_0}u\Big\}\Big\rceil,\\
G=d\sigma K_0,\qquad K_0=\max\{1,\tfrac M\sigma,\sqrt V\},\qquad \Gamma=M\big(1+\tfrac Rr\big).
\end{gathered}
\end{equation}

\paragraph{Configurations and certificates.} For a block with $n$ updates, the two
configurations are (according to \ref{eq:probit-configurations2})
\begin{equation}\label{eq:probit-configurations}
\begin{aligned}
\mathrm c:&\quad h=\min\Big\{\frac r2,\frac{RG}{\Gamma\sqrt n}\Big\},\quad \eta_t=\frac{R}{G\sqrt n},\quad w_t=1;\\
\mathrm s:&\quad h=\min\Big\{\frac r2,\frac{G^2}{\mu\Gamma(n+1)}\Big\},\quad \eta_t=\frac{2}{\mu(t+1)},\quad w_t=t,
\end{aligned}
\end{equation}
and their certificates, i.e.\ computable upper bounds on the error of the block output, are
\begin{equation}\label{eq:probit-selection}
U_0(b)=MR,\qquad U_{\mathrm c}(b)=\frac{2RG}{\sqrt n}+MRe^{-ub/2},\qquad
U_{\mathrm s}(b)=\frac{3G^2}{\mu(n+1)}+MRe^{-ub/2}.
\end{equation}
\begin{algorithm}[ht]
\caption{Debiased probit comparisons with doubling}
\label{alg:probit}
\begin{algorithmic}[1]
\Require anchor $c=0$ with $B(0,r)\subseteq Q\subseteq B(0,R)$; unit vector $e_0$; constants \eqref{eq:probit-constants}; coefficients $a_j$ of \cref{lem:probit-taylor}
\Statex \textbf{Output:} $\hat x$ at the moment the budget of calls is exhausted
\State $\hat x\gets0$
\For{$\ell=0,1,2,\dots$}
  \State $b\gets2^\ell$ \Comment{calls in this block}
  \State $n\gets\lfloor b/C_0\rfloor$ \Comment{number of updates}
  \State $\tau\gets0$, \ $z\gets0$ \Comment{calls used so far, output of the block}
  \State $\mathrm{cfg}\gets\arg\min\{U_0(b),U_{\mathrm c}(b),U_{\mathrm s}(b)\}$ \Comment{\eqref{eq:probit-selection}; $U_{\mathrm c}=U_{\mathrm s}=\infty$ if $n=0$, $U_{\mathrm s}=\infty$ if $\mu=0$}
  \If{$\mathrm{cfg}\in\{\mathrm c,\mathrm s\}$}
    \State $(h,\eta_t,w_t)\gets$ configuration $\mathrm{cfg}$ of \eqref{eq:probit-configurations}, \ $x_1\gets0$, \ $\mathrm{done}\gets\mathrm{true}$
    \For{$t=1,\dots,n$}
      \State draw $e_t\sim\Unif(\S^{d-1})$ and $J_t$ with $\Pr(J_t=j)=p_j$
      \If{$\tau+kJ_t>b$}
        \State $\mathrm{done}\gets\mathrm{false}$; \textbf{break} \Comment{number of the block exceeded}
      \EndIf
      \For{$i=1,\dots,J_t$}
        \State $B_{t,i}\gets\frac1k\sum_{j=1}^kO(x_t+he_t,\,x_t-he_t)$ \Comment{$k$ fresh calls}
      \EndFor
      \State $T_t\gets\frac{a_{J_t}}{p_{J_t}}\prod_{i=1}^{J_t}B_{t,i}$, \ $\tau\gets\tau+kJ_t$ \Comment{unbiased estimate, \cref{lem:probit-estimator}}
      \State $x_{t+1}\gets\Projp_{K_h}\big(x_t-\eta_t\,d\sigma\,T_te_t\big)$ \Comment{$K_h=(1-h/r)Q$}
    \EndFor
    \If{$\mathrm{done}$}
      \State $z\gets\sum_{t=1}^nw_tx_t\big/\sum_{t=1}^nw_t$
    \EndIf
  \EndIf
  \For{$i=1,\dots,b-\tau$}
    \State $O(0,\frac r2e_0)$ \Comment{dummy call, answer ignored}
  \EndFor
  \State $\hat x\gets z$
\EndFor
\end{algorithmic}
\end{algorithm}

\section{Linear link for Preference Oracle}
\subsection{Algorithm pseudo-code}
\begin{algorithm}[t]
\caption{Anchored value reduction with doubling}
\label{alg:value}
\begin{algorithmic}[1]
\Require $c\in Q$, $r>0$ with $B(c,r)\subseteq Q$; $\mathsf{BCO}$ (Algorithm 1 of \cite{fokkema2024online})
\Statex \textbf{Output:} $\hat x$ at the moment the budget of calls is exhausted
\State $\hat x\gets c$
\For{$\ell=0,1,2,\dots$}
  \State $b\gets2^\ell$ \Comment{horizon}
  \State $\delta\gets(b+1)^{-2}$ \Comment{confidence level}
  \State draw a fresh seed $U_\ell$ \Comment{internal randomness of $\mathsf{BCO}$}
  \For{$t=1,\dots,b$}
    \State $x_t\gets\mathsf{BCO}(Q,b,\delta;\,U_\ell;\,y_1,\dots,y_{t-1})$ \Comment{point played in round $t$ after observing losses $y_1,\dots,y_{t-1}$}
    \If{$x_t\ne c$}
      \State $o_t\gets O(x_t,c)$, \ $y_t\gets\frac{1+o_t}2$
    \Else
      \State $O(c+\frac r2e_1,c)$ \Comment{dummy call, answer ignored}
      \State $y_t\sim\Unif\{0,1\}$
    \EndIf
  \EndFor
  \State $\hat x\gets\frac1b\sum_{t=1}^bx_t$
\EndFor
\end{algorithmic}
\end{algorithm}

Algorithm~1 of \citet{fokkema2024online} is an interactive procedure: it is initialized with
the set $Q$, a horizon $n$ and a confidence level $\delta$, and then in every round it plays a
point and receives a noisy value of the loss at that point. All its remaining tuning
parameters are explicit functions of $n$, $d$, $\delta$ and the geometry of $Q$ (in the
stochastic setting, \cite[Section~5]{fokkema2024online}), and the affine change of
coordinates it applies to $Q$ is internal to the routine. We write
$\mathsf{BCO}(Q,b,\delta;U_\ell;y_1,\dots,y_{t-1})$ for the point it plays in round $t$ when
its internal randomness is $U_\ell$ and the values received in rounds $1,\dots,t-1$ are
$y_1,\dots,y_{t-1}$. Computationally, this simply continues the same run. In their notation,
$x_t$ is the played point $A_t$, and $y_t$ is the observed loss $\omega(A_t)+\varepsilon_t$ for the given objective function $\omega$ and noise $\varepsilon$ (in the proof below, $\omega = g$ and $\varepsilon=\zeta$). The transformed observation $Y_t$ of
their Algorithm~1 is computed by the routine itself from $y_t$. Each block uses a fresh seed,
so the blocks are independent and \cref{lem:value} applies to each of them separately. The
oracle is defined only for distinct points, so if $\mathsf{BCO}$ requests $x_t=c$, we make a
dummy call at $(c+\frac r2e_1,c)$, which lies in $Q\times Q$ because $B(c,r)\subseteq Q$
($e_1$ may be any unit vector), and feed a fair coin, whose mean $\frac12=g(c)$ matches the
other rounds. This keeps the number of calls of block $\ell$ equal to $2^\ell$, so that the
last completed block after $N$ calls is determined by $N$ alone. Because of the doubling, the
method needs neither $N$ nor $\mu$, $\lambda$, $M$, $R$: at every moment it holds the average
point of the last completed block and returns it when the budget is exhausted.

\subsection{Proof of \Cref{thm:lin-upper}}\label{proof-thm-lin-upper}
Assume that $B(c,r)\subseteq Q\subseteq B(c,R)$ for some $0<r\le R$. Let $x^*\in\arg\min_{x\in Q}f$ and $f^*=f(x^*)$. Since $f$ is $M$-Lipschitz on $Q$ and $c,x^*\in Q\subseteq B(c,R)$,
\begin{equation}\label{eq:lin-trivial}
  f(c)-f^*\le MR,\qquad f(x)-f^*\le 2MR\quad\text{for all }x\in Q .
\end{equation}
Let $C_0,p>0$ be the constants of \cref{lem:value}. Since $p>0$,
$s\mapsto\log^p s$ is nondecreasing on $[e,\infty)$. Every iteration of the inner loop of
\cref{alg:value} makes exactly one call: if $x_t\ne c$, it queries $(x_t,c)$, and if $x_t=c$,
it makes a dummy call at $(c+\frac r2e_1,c)$, which is feasible since $B(c,r)\subseteq Q$, and
ignores its answer. Hence, the $\ell$-th iteration of the outer loop (the $\ell$-th block) uses
exactly $b=2^\ell$ calls. For one iteration of the outer loop, define
\begin{equation}\label{eq:lin-cert}
  U_b:=\frac{2C_0\,d^{3/2}\log^{p}\!\big(ed(b+1)^3\big)}{\lambda\sqrt b}+\frac{2MR}{(b+1)^2},
\end{equation}
and let the block publish $z_b=\frac1b\sum_{t\le b}x_t$ if $U_b<MR$, and $c$ otherwise.

Fix a block and an iteration $t$ of its inner loop. Denote by $\cH_t$ the history before this
iteration: the internal randomness
of $\mathsf{BCO}$ in iterations $1,\dots,t-1$, and the points $x_s$, the oracle answers and the
observations $y_s$ of iterations $s<t$.

\medskip
\emph{Step 1.} Define $g:Q\to\R$ by
\[
  g(x):=\frac12+\frac\lambda2\big(f(x)-f(c)\big).
\]
Since $\lambda>0$, $g$ is convex, $g(c)=\frac12$, $\arg\min_Q g=\arg\min_Qf$, and
\begin{equation}\label{eq:lin-ell-gap}
  g(x)-\min_Q g=\frac\lambda2\big(f(x)-f^*\big).
\end{equation}
Moreover, $|f(x)-f(c)|\le M\|x-c\|\le MR$ on $Q$, hence $|g(x)-\frac12|\le\frac{\lambda MR}2\le\frac14$, i.e.
$g(Q)\subseteq[\frac14,\frac34]\subseteq[0,1]$.

Let $x_t\in Q$ be the point requested by $\mathsf{BCO}$ at iteration $t$. If $x_t\ne c$,
the algorithm queries $(x_t,c)$ and feeds $y_t=\frac{1+o_t}2\in\{0,1\}$. Since the oracle uses
fresh randomness at every call and $|f(x_t)-f(c)|\le MR$, \cref{saha_oracle_nonlinear} and
\eqref{eq:linear-link} give
\[
  \E[y_t\mid\cH_t,x_t]=\frac12+\frac\lambda2\big(f(x_t)-f(c)\big)=g(x_t).
\]
If $x_t=c$, then $y_t$ is a fair coin independent of $\cH_t$ and of the dummy answer, and
$\E[y_t\mid\cH_t,x_t]=\frac12=g(c)$.
So, in both cases $y_t=g(x_t)+\zeta_t$ with $\E[\zeta_t\mid\cH_t,x_t]=0$. Since $y_t\in\{0,1\}$ and
$g(x_t)\in[\frac14,\frac34]$, we have $|\zeta_t|\le\frac34$. Therefore,
\[
  \E\big[e^{\zeta_t^2}\,\big|\,\cH_t,x_t\big]\le e^{9/16}<2 .
\]
Thus, the observations satisfy the assumptions of \cref{lem:value} with the convex function $g:Q\to[0,1]$.

\medskip
\emph{Step 2.} Fix a block with budget $b$, run with confidence $\delta:=(b+1)^{-2}$. By \cref{lem:value} (with $n=b$) and
\eqref{eq:lin-ell-gap}, there is an event $\mathcal G$ with $\Pr(\mathcal G)\ge1-\delta$ on which
\[
  \sum_{t\le b}\big(f(x_t)-f^*\big)=\frac2\lambda\sum_{t\le b}\big(g(x_t)-\min_Q g\big)
  \le\frac{2C_0\,d^{3/2}\sqrt b}{\lambda}\log^p\!\Big(\frac{edb}{\delta}\Big)
  \le\frac{2C_0\,d^{3/2}\sqrt b}{\lambda}\log^p\!\big(ed(b+1)^3\big),
\]
where we used $b/\delta=b(b+1)^2\le(b+1)^3$. The point $z_b=\frac1b\sum_{t\le b}x_t$ lies in $Q$ by convexity, and
Jensen's inequality gives, on $\mathcal G$,
\[
  f(z_b)-f^*\le\frac1b\sum_{t\le b}\big(f(x_t)-f^*\big)\le\frac{2C_0\,d^{3/2}\log^p(ed(b+1)^3)}{\lambda\sqrt b}.
\]
On $\mathcal G^c$ we use $f(z_b)-f^*\le2MR$ from \eqref{eq:lin-trivial}. Therefore
\begin{equation}\label{eq:lin-block-z}
  \E f(z_b)-f^*\le\frac{2C_0\,d^{3/2}\log^p(ed(b+1)^3)}{\lambda\sqrt b}+2MR\,\delta=U_b .
\end{equation}

\medskip
\emph{Step 3.}
Let $\hat z_b$ be the point published by the block. If $U_b<MR$, then $\hat z_b=z_b$ and
$\E f(\hat z_b)-f^*\le U_b$ by \eqref{eq:lin-block-z}. If $U_b\ge MR$, then $\hat z_b=c$ and
$f(c)-f^*\le MR$ by \eqref{eq:lin-trivial}. In both cases,
\begin{equation}\label{eq:lin-block}
  \E f(\hat z_b)-f^*\le\min\{MR,\,U_b\}.
\end{equation}

\medskip
\emph{Step 4.}
Fix $N\ge1$ and let $\ell:=\lfloor\log_2(N+1)\rfloor-1\ge0$. Blocks $0,\dots,\ell$ use $2^{\ell+1}-1\le N$ calls, while
blocks $0,\dots,\ell+1$ would need $2^{\ell+2}-1>N$ calls. Hence after $N$ calls the published point $\hat x_N$ is the output of
block $\ell$, whose budget satisfies
\[
  b=2^\ell\in\Big(\frac{N+1}4,\frac{N+1}2\Big].
\]
By \eqref{eq:lin-block}, $\E f(\hat x_N)-f^*\le\min\{MR,U_b\}$, and it remains to bound $U_b$.

Since $\lambda MR\le\frac12$ we have $\frac{2MR}{(b+1)^2}\le\frac1{\lambda(b+1)^2}\le\frac1{\lambda\sqrt b}$.
Since $b\le\frac{N+1}2$ and $N\ge1$, we get $b+1\le N+1$, and as $ed\ge1$,
$\log\big(ed(b+1)^3\big)\le\log\big((ed(N+1))^3\big)=3\log\big(ed(N+1)\big)$. Finally $b>\frac{N+1}4$ gives
$\frac1{\sqrt b}<\frac2{\sqrt N}$. Using also $d\ge1$ and $\log\big(ed(N+1)\big)\ge1$,
\[
  U_b\le\frac{d^{3/2}}{\lambda\sqrt b}\Big(2C_03^p\log^p\!\big(ed(N+1)\big)+1\Big)
  \le\frac{2\,(2C_03^p+1)\,d^{3/2}\log^p\!\big(ed(N+1)\big)}{\lambda\sqrt N}.
\]
Setting $C:=2(2C_03^p+1)$, we conclude
\[
  \E f(\hat x_N)-f^*\le\min\Big\{MR,\ \frac{C\,d^{3/2}\log^p\big(ed(N+1)\big)}{\lambda\sqrt N}\Big\}
\]

Notice that we did not assume strong convexity, so it holds for any $\mu\ge 0$.
\providecommand{\KL}{\operatorname{KL}}
\providecommand{\TV}{\operatorname{TV}}

\providecommand{\KL}{\operatorname{KL}}
\providecommand{\TV}{\operatorname{TV}}

\section{Lower bounds for linear link}
\label{proof:thm-lin-lower}
All subsections in this section share common definitions.
\subsection{Preliminaries}
Let us recall the notation for an oracle $o$. An algorithm $A\in\cA_N$ consists of a
probability space $(\mathcal U,\nu)$ (its internal randomness) and measurable maps
$\phi_1,\dots,\phi_{N+1}$. Run on an oracle $o^f$ for an objective $f$ (we write $o$ when $f$
is clear), it draws a seed $U\sim\nu$ independent of the oracle and, for $t=1,\dots,N$,
queries the pair
\[
  (x_t,y_t):=\phi_t(U,o_1,\dots,o_{t-1})\in Q\times Q
\]
and receives the $t$-th answer $o_t:=o^f(x_t,y_t)$. Its output is
$\hat x:=\phi_{N+1}(U,o_1,\dots,o_N)\in Q$. An algorithm that makes fewer than $N$ calls is
included by padding it with calls whose answers it ignores.

The \emph{transcript} of $A$ is $T:=(U,o_1,\dots,o_N)$; note that $x_t$, $y_t$ and $\hat x$ are
functions of $T$. For a query $(x,y)$, let $m_f(x,y):=\lambda\big(f(x)-f(y)\big)$ be the mean
answer of $o^f$, and let
\[
  q_f(s\mid x,y):=\frac{1+s\,m_f(x,y)}{2},\qquad s\in\{-1,+1\},
\]
be the probability that the answer equals $s$, i.e.\ $q_f(\pm1\mid x,y)=\frac{1\pm m_f(x,y)}2$.
For a fixed algorithm $A$, $P_f$ denotes the law of $T$ when $A$ is run on $o^f$: $U\sim\nu$
and, given $U$, the answers $(o_1,\dots,o_N)$ have probability
$\prod_{t=1}^N q_f(o_t\mid x_t,y_t)$. $\E_f$ denotes the expectation over $T\sim P_f$.
\subsection{Auxiliary lemmas}
\begin{lemma}\label{lem:lin-kl}
Fix an algorithm $A\in\cA_N$ and $\varpi\in[0,1)$. Let $f,g:Q\to\R$ satisfy
$\lambda|f(x)-f(y)|\le1$ and $\lambda|g(x)-g(y)|\le\varpi$ for all $x,y\in Q$, and put $\Delta_{f,g}:=f-g$. Then,
\[
  \KL(P_f\,\|\,P_g)\ \le\ \frac{\lambda^2}{1-\varpi^2}\;\E_f\sum_{t=1}^N\big(\Delta_{f,g}(x_t)-\Delta_{f,g}(y_t)\big)^2.
\]

\end{lemma}

\begin{proof}
Under $P_f$ and $P_g$ the seed has the same law, and given the past
$(U,o_1,\dots,o_{t-1})$ the query $(x_t,y_t)$ is the same under both. Hence, the chain
rule for the KL divergence gives
\[
  \KL(P_f\|P_g)=\E_f\sum_{t\le N}\KL\big(q_f(\cdot\mid x_t,y_t)\,\|\,q_g(\cdot\mid x_t,y_t)\big).
\]
Fix $t$ and write $p:=\lambda(f(x_t)-f(y_t))$, $p':=\lambda(g(x_t)-g(y_t))$ for the two
means. For distributions $q,q'$ on $\{\pm1\}$ with means $p,p'$ and $|p'|<1$, Jensen's
inequality gives
\[
  \KL(q\|q')\le\log\E_q\frac{q}{q'}=\log\big(1+\chi^2(q\|q')\big)\le\chi^2(q\|q'),
\]
and the chi-squared divergence is explicit:
\[
  \chi^2(q\|q')
  =\sum_{s=\pm1}\frac{(q(s)-q'(s))^2}{q'(s)}
  =\frac{(p-p')^2}{4}\Big(\frac{2}{1+p'}+\frac{2}{1-p'}\Big)
  =\frac{(p-p')^2}{1-p'^2}.
\]
By definition of $\Delta_{f,g}$, $p-p'=\lambda(\Delta_{f,g}(x_t)-\Delta_{f,g}(y_t))$, while $|p'|\le\varpi$. Substituting
into the bound and summing over $t$ yields the claim.
\end{proof}

\begin{lemma}\label{lem:assouad}
Fix an algorithm $A\in\cA_N$. Let $\{f_v\}_{v\in\{\pm1\}^d}$ be functions $Q\to\R$ such that
$\lambda|f_v(x)-f_v(y)|\le1$ for all $x,y\in Q$, and let $P_v$ be the law of the
transcript of $A$ run on $o^{f_v}$, with $\E_v$ the expectation under $P_v$.
Suppose $\ell:\R\times\{\pm1\}\to[0,\infty)$ and $\Delta>0$ satisfy, for every $v$,
\begin{enumerate}
  \item[\textup{(i)}] $f_v(x)-\min_Qf_v=\sum_{i=1}^d\ell(x_i,v_i)$ for all $x\in Q$;
  \item[\textup{(ii)}] $\ell(t,+1)+\ell(t,-1)\ge\Delta$ for all $t\in\R$.
\end{enumerate}
Let $v^{(i)}$ be $v$ with the $i$-th coordinate flipped and
\[
  \bar K:=\frac1{d\,2^d}\sum_{v\in\{\pm1\}^d}\sum_{i=1}^d\KL\big(P_v\,\|\,P_{v^{(i)}}\big).
\]
Then,
\[
  \max_v\E_v\big[f_v(\hat x)-\min_Qf_v\big]\ \ge\ \frac{d\Delta}{2}\Big(1-\sqrt{\bar K/2}\Big).
\]
\end{lemma}
\begin{proof}
By (i), and since the maximum is at least the average,
\[
  \max_v\E_v\big[f_v(\hat x)-\min_Qf_v\big]
  \;\ge\;\sum_i 2^{-d}\sum_v\E_v\,\ell(\hat x_i,v_i).
\]
Fix $i$. For $v$ with $v_i=+1$, using $\ell\ge0$ and (ii),
\[
  \E_v\,\ell(\hat x_i,+1)+\E_{v^{(i)}}\ell(\hat x_i,-1)
  \;\ge\;\int\big[\ell(\hat x_i,+1)+\ell(\hat x_i,-1)\big]\,d\big(P_v\wedge P_{v^{(i)}}\big)
  \;\ge\;\Delta\big(1-\TV(P_v,P_{v^{(i)}})\big),
\]
where $P_v\wedge P_{v^{(i)}}$ is the measure with density
$\min\{dP_v/d\rho,\,dP_{v^{(i)}}/d\rho\}$ with respect to a dominating measure $\rho$, so
its total mass is $1-\TV(P_v,P_{v^{(i)}})$. Since $\hat x$ is a function of the
transcript, and since $v\mapsto v^{(i)}$ is a bijection from $\{v_i=+1\}$ to $\{v_i=-1\}$
while $\TV$ is symmetric,
\[
  2^{-d}\sum_v\E_v\,\ell(\hat x_i,v_i)\;\ge\;\frac\Delta2\,(1-T_i),
  \qquad
  T_i:=2^{-d}\sum_v\TV(P_v,P_{v^{(i)}}).
\]
Summing over $i$ gives $\frac{d\Delta}2\big(1-\frac1d\sum_iT_i\big)$. By Pinsker's
inequality, $\TV(P_v,P_{v^{(i)}})\le\sqrt{\KL(P_v\|P_{v^{(i)}})/2}$, and since
$\sqrt{\cdot}$ is concave, Jensen's inequality gives
$\frac1d\sum_iT_i\le\sqrt{\bar K/2}$.
\end{proof}

\subsection{Proof of \cref{thm:lin-lower}.}

\paragraph{Case $\mu=0$.}
Since the class $\Fzero$ and the oracle are invariant under translations, we may assume
$c=0$, so that $\B^d(r)\subseteq Q\subseteq\B^d(R)$. Fix $N\ge1$ and an algorithm
$A\in\cA_N$, and put
\[
  a:=\frac{M}{2\sqrt d},\qquad
  \tau:=\min\Big\{\frac{Mr}{4d},\ \frac{\sqrt3}{16\lambda\sqrt N}\Big\},\qquad
  s:=\frac{\tau}{a}.
\]
For $v\in\{\pm1\}^d$ define
\[
  f_v(x):=\sum_{i=1}^d\ell(x_i,v_i),\qquad \ell(t,w):=a|t-ws|.
\]

Each $\ell(\cdot,w)$ is convex, so $f_v$ is convex. Every subgradient of
$f_v$ has coordinates in $[-a,a]$, hence norm at most $a\sqrt d=M/2$, and $f_v\in\Fzero$. Since
$\ell(\cdot,w)\ge0=\ell(ws,w)$, the point $sv$ minimizes $f_v$ and $\min_Qf_v=0$. As
$\|sv\|=s\sqrt d=2d\tau/M\le r/2$, it lies in $Q$. Finally, for $x,y\in Q$,
$\lambda|f_v(x)-f_v(y)|\le\lambda\frac M2\|x-y\|\le\lambda MR\le\frac12$, so all gaps lie in the
linear range of \eqref{eq:linear-link}, and the hypotheses of \cref{lem:lin-kl,lem:assouad}
hold with $\varpi=\frac12$.

Condition (i) of \cref{lem:assouad} holds by construction, and by the
triangle inequality $\ell(t,+1)+\ell(t,-1)=a\big(|t-s|+|t+s|\big)\ge2as$ for all $t$, so (ii)
holds with $\Delta=2as=2\tau$.

The functions $f_v$ and $f_{v^{(i)}}$ differ only in the $i$-th term, so
\[
  \Delta_{f_v,f_{v^{(i)}}}(x)=a\big(|x_i-v_is|-|x_i+v_is|\big)=-2av_i\max\{-s,\min\{s,x_i\}\},
\]
and hence $|\Delta_{f_v,f_{v^{(i)}}}(x)-\Delta_{f_v,f_{v^{(i)}}}(y)|\le4as=4\tau$ for all
$x,y\in Q$. By \cref{lem:lin-kl} with $\kappa:=1/(1-\varpi^2)=\frac43$,
$\KL(P_v\|P_{v^{(i)}})\le16\kappa N\lambda^2\tau^2$ for all $v,i$, so
$\bar K\le16\kappa N\lambda^2\tau^2=\frac{64}3N\lambda^2\tau^2$. Since
$\tau\le\frac{\sqrt3}{16\lambda\sqrt N}$, we get
$\sqrt{\bar K/2}\le\lambda\tau\sqrt{32N/3}\le\frac{\sqrt2}4\le\frac12$.

By \cref{lem:assouad},
\[
  \max_v\E_v\big[f_v(\hat x)-\min_Qf_v\big]
  \ge\frac{d\Delta}2\Big(1-\sqrt{\bar K/2}\Big)
  \ge\frac{d\tau}2
  =\min\Big\{\frac{Mr}8,\ \frac{\sqrt3\,d}{32\lambda\sqrt N}\Big\}.
\]
Since $A\in\cA_N$ was arbitrary and $f_v\in\Fzero$ for every $v$, the claim follows.

\medskip\noindent\textbf{Case $0<\mu\le M/(8R)$.}
Let $\delta>0$ be defined by
\[
  \delta^2:=\min\Big\{\frac{r^2}{3d},\ \frac{d}{36\,\lambda^2\mu^2R^2N}\Big\},
\]
and for $v\in\{\pm1\}^d$ define
\[
  f_v(x):=\sum_{i=1}^d\ell(x_i,v_i)=\frac\mu2\|x-\delta v\|^2,\qquad
  \ell(t,w):=\frac\mu2(t-\delta w)^2.
\]

The function $f_v$ is $\mu$-strongly convex. Since
$\|\delta v\|=\delta\sqrt d\le r/\sqrt3\le r$, the point $\delta v$ lies in $Q$, so
$\min_Qf_v=0$. On $Q$, $\|\nabla f_v(x)\|=\mu\|x-\delta v\|\le2\mu R\le M/4$, so $f_v\in\Fm$.
Moreover, $0\le f_v\le2\mu R^2\le MR/4$ on $Q$, so $\lambda|f_v(x)-f_v(y)|\le\lambda MR/4\le\frac18$;
hence all gaps lie in the linear range of \eqref{eq:linear-link}, and the hypotheses of
\cref{lem:lin-kl,lem:assouad} hold with $\varpi=\frac18$.

Condition (i) of \cref{lem:assouad} holds by construction, and
$\ell(t,+1)+\ell(t,-1)=\mu(t^2+\delta^2)\ge\mu\delta^2$ for all $t$, so (ii) holds with
$\Delta=\mu\delta^2$.

 Since $\|v\|=\|v^{(i)}\|$,
\[
  \Delta_{f_v,f_{v^{(i)}}}(x)=-2\mu\delta\,v_ix_i,
\]
and hence
$\sum_i\big(\Delta_{f_v,f_{v^{(i)}}}(x)-\Delta_{f_v,f_{v^{(i)}}}(y)\big)^2=4\mu^2\delta^2\|x-y\|^2\le16\mu^2\delta^2R^2$
for all $x,y\in Q$. Applying \cref{lem:lin-kl} with $\kappa':=1/(1-\varpi^2)=64/63\le9/8$ for
each $i$ and summing over $i$, we get
$\sum_i\KL(P_v\|P_{v^{(i)}})\le16\kappa'\lambda^2\mu^2R^2N\delta^2$ for every $v$, so
$\bar K\le16\kappa'\lambda^2\mu^2R^2N\delta^2/d$. Since $\delta^2\le\frac{d}{36\lambda^2\mu^2R^2N}$,
we get $\bar K\le4\kappa'/9\le\frac12$, hence $\sqrt{\bar K/2}\le\frac12$.

By \cref{lem:assouad},
\[
  \max_v\E_v\big[f_v(\hat x)-\min_Qf_v\big]
  \ge\frac{d\Delta}2\Big(1-\sqrt{\bar K/2}\Big)
  \ge\frac{d\mu\delta^2}{4}
  =\min\Big\{\frac{\mu r^2}{12},\ \frac{d^2}{144\,\lambda^2\mu R^2N}\Big\}.
\]

\medskip
Combining the two bounds for $0\le\mu\le M/(8R)$, we get the final bound
\[
  \cL_N(\Fm,o)\ \ge\ \max\Big\{
    \min\Big\{\frac{21Mr}{256},\frac{3\sqrt3d}{128\lambda\sqrt N}\Big\},\
    \min\Big\{\frac{\mu r^2}{12},\frac{d^2}{144\lambda^2\mu R^2N}\Big\}\Big\},
\]
which gives the claimed lower bound.

\subsection{Proof of theorem \ref{thm:restricted}}\label{proof:restricted}
 The proofs of
\cref{lem:lin-kl,lem:assouad} use only that the queries and the output are functions of the
transcript and the bounds on $\lambda|f(x)-f(y)|$ at the queried pairs, so both lemmas apply
to this class of algorithms \ref{thm:restricted}. We use \cref{lem:assouad} with (i) in the weaker form
$f_v(x)-\min_Qf_v\ge\sum_i\ell(x_i,v_i)$, which is all its proof uses.

For $\mu\ge0$, $b>0$ and $v\in\{\pm1\}^d$ let
\[
  f_v(x):=\frac\mu2\|x\|^2-b\ip vx,
\]
and assume that
\begin{equation}\label{eq:restricted-range}
  \lambda H\big(\mu R+b\sqrt d\big)\le\tfrac14 .
\end{equation}
At the $k$-th query the mean answer is
\[
  \lambda\big(f_v(x_k+h_ke_k)-f_v(x_k-h_ke_k)\big)=2\lambda h_k\big(\mu\ip{x_k}{e_k}-b\ip v{e_k}\big),
\]
and its absolute value is at most $2\lambda H(\mu R+b\sqrt d)\le\frac12$, since $\|x_k\|\le R$,
$|\ip v{e_k}|\le\sqrt d$ and \eqref{eq:restricted-range} holds. Thus, the hypotheses of
\cref{lem:lin-kl,lem:assouad} hold at all queried pairs with $\varpi=\frac12$. The quadratic
part cancels between neighbors: $\Delta_{f_v,f_{v^{(j)}}}(x)=-2bv_jx_j$, so its increment over
the $k$-th pair is $-4bv_jh_ke_{k,j}$. By \cref{lem:lin-kl} with $\kappa=1/(1-\varpi^2)=\frac43$
and $\sum_je_{k,j}^2=1$,
\[
  \sum_{j=1}^d\KL(P_v\|P_{v^{(j)}})
  \le\frac43\lambda^2\,\E_v\sum_{k=1}^N\sum_{j=1}^d16b^2h_k^2e_{k,j}^2
  =\frac{64}3\lambda^2b^2\,\E_v\sum_{k=1}^Nh_k^2
  \le\frac{64}3\lambda^2b^2B_N
\]
for every $v$, and hence $\bar K\le\frac{64}{3d}\lambda^2b^2B_N$.

\paragraph{Case $\mu=0$.} Take $\mu=0$ and
$b:=\min\{\frac M{2\sqrt d},\frac{\sqrt d}{16\lambda\sqrt{B_N}}\}$. Then
\eqref{eq:restricted-range} holds, since $\lambda Hb\sqrt d=\frac{b\sqrt d}{2M}\le\frac14$. Each
$f_v$ is linear with Lipschitz constant $b\sqrt d\le\frac M2$, so $f_v\in\Fzero$, and it is
minimized over $Q$ at $x_v^*:=Rv/\sqrt d$. For $x\in Q$, since $\|x\|\le R$,
\[
  f_v(x)-\min_{x\in Q}f_v=b\sqrt d\Big(R-\Big\langle\frac v{\sqrt d},x\Big\rangle\Big)
  \ge\frac{b\sqrt d}{2R}\|x-x_v^*\|^2 .
\]
If $v_jx_j\le0$, then $x_j$ and $Rv_j/\sqrt d$ have opposite signs or $x_j=0$, so
$|x_j-Rv_j/\sqrt d|\ge R/\sqrt d$. Hence, for $x\in Q$,
\[
  f_v(x)-\min_Qf_v\ge\frac{b\sqrt d}{2R}\sum_{j=1}^d\Big(x_j-\frac{Rv_j}{\sqrt d}\Big)^2
  \ge\frac{b\sqrt d}{2R}\cdot\frac{R^2}{d}\sum_{j=1}^d\ind\{v_jx_j\le0\}
  =\sum_{j=1}^d\ell(x_j,v_j),
\]
where $\ell(t,w):=\frac{bR}{2\sqrt d}\ind\{wt\le0\}$, so (i) holds. Moreover,
$\ell(t,+1)+\ell(t,-1)=\frac{bR}{2\sqrt d}\big(\ind\{t\le0\}+\ind\{t\ge0\}\big)\ge\frac{bR}{2\sqrt d}$
for every $t$, so (ii) holds with $\Delta=\frac{bR}{2\sqrt d}$. Since
$b\le\frac{\sqrt d}{16\lambda\sqrt{B_N}}$, we have $b^2\le\frac{d}{256\lambda^2B_N}$, and the
bound on $\bar K$ above gives
$\bar K\le\frac{64}{3d}\lambda^2B_N\cdot\frac{d}{256\lambda^2B_N}=\frac1{12}$, so
$\sqrt{\bar K/2}\le\frac1{\sqrt{24}}\le\frac12$, and \cref{lem:assouad} yields
\[
  \max_v\E_v\big[f_v(\hat x)-\min_Qf_v\big]\ge\frac{d\Delta}4=\frac{bR\sqrt d}8
  =\min\Big\{\frac{MR}{16},\ \frac{dR}{128\lambda\sqrt{B_N}}\Big\}.
\]

\paragraph{Case $\mu>0$.}  Let $\bar M:=\min\{\mu R,\frac{M-2\mu R}2\}$ and take
$b:=\min\{\frac{\bar M}{\sqrt d},\frac{\sqrt{3d}}{8\sqrt2\,\lambda\sqrt{B_N}}\}$. Then
$\mu R+b\sqrt d\le\mu R+\bar M\le\frac M2$, so \eqref{eq:restricted-range} holds. On $\B^d(2R)$,
$\|\nabla f_v(x)\|=\|\mu x-bv\|\le2\mu R+\bar M\le M$, so $f_v\in\Fm(\B^d(2R))$. Since
$b\sqrt d\le\mu R$, the minimizer $(b/\mu)v$ lies in $Q$, and
\[
  f_v(x)-\min_{x\in Q}f_v=\frac\mu2\Big\|x-\frac b\mu v\Big\|^2 .
\]
If $v_jx_j\le0$, then $|x_j-(b/\mu)v_j|\ge b/\mu$; hence (i) holds with
$\ell(t,w):=\frac{b^2}{2\mu}\ind\{wt\le0\}$, and (ii) holds with $\Delta=\frac{b^2}{2\mu}$. The
choice of $b$ gives $\bar K\le\frac12$, so $\sqrt{\bar K/2}\le\frac12$, and \cref{lem:assouad}
yields
\[
  \max_v\E_v\big[f_v(\hat x)-\min_Qf_v\big]\ge\frac{d\Delta}4=\frac{db^2}{8\mu}
  =\min\Big\{\frac{\bar M^2}{8\mu},\ \frac{3d^2}{1024\lambda^2\mu B_N}\Big\}.
\]

\begin{remark}
    If $h_{\max}\le H$, then \cref{alg:ardg} with $N$
iterations and $m$ repetitions belongs to this class with $mN$ calls and radius budget
$m\sum_{k<N}h_k^2$. If $h_k\le c(k+1)^{-1/4}$, then by \cref{lem:sums}(a) this budget is at
most $2mc^2\sqrt N$, and the convex bound becomes
$\min\{\frac{MR}{16},\frac{dR}{128\sqrt2\lambda c\sqrt m}N^{-1/4}\}$; if
$h_k\le c(k+1)^{-1/3}$, it is at most $3mc^2N^{1/3}$, and the strongly convex bound becomes
$\min\{\frac{\bar M^2}{8\mu},\frac{d^2}{1024\lambda^2\mu mc^2}N^{-1/3}\}$.
\end{remark} 

\section{Upper bounds for general link}
In this section, all the lemmas and theorems share common definitions, according to \cref{alg:ardg}, \cref{thm:ardg-convex}, and \cref{thm:ardg-sc}.
\begin{lemma}
    \label{lem:estimator}
Let $f\in \Fzero(U)$. Then,

\begin{enumerate}
  \item[(i)] $\E[g_k\mid\Fk]=a_k\nabla f_{h_k}(x_k)+b_k$, where
  \[
  a_k:=\frac{2\lambda h_k}{d},\qquad
    b_k:=\E\big[r(\Delta_k)e_k\mid\Fk\big],
    \qquad
    \norm{b_k}\le\frac{\rho(2Mh_k)}{\sqrt d}.
  \]

  \item[(ii)] $\norm{g_k}\le1$.

  \item[(iii)] $\E\big[\norm{g_k}^2\mid\Fk,e_k\big]
  \le G_k^2:=\frac1m+\Big(1-\frac1m\Big)\Psi(2Mh_k)^2$.
\end{enumerate}
\end{lemma}
\begin{proof}
Since $f$ is $M$-Lipschitz and the two query points are at distance $2h_k$,
\[
  |f(x_k+h_ke)-f(x_k-h_ke)|\le 2Mh_k,
\]
Further, we will denote $\Delta_k := f(x_k+h_ke)-f(x_k-h_ke)$.
\emph{(i).} Conditionally on $\Fk$ and $e_k$, the $m$ answers are i.i.d.\ with mean
$p_k=\psi(\Delta_k)$, so $\E[\bar O_k\mid\Fk,e_k]=\psi(\Delta_k)$. By the tower property,
\[
  \E[g_k\mid\Fk]=\E[\psi(\Delta_k)e_k\mid\Fk].
\]
Writing $\psi(t)=\lambda t+r(t)$,
\[
  \E[g_k\mid\Fk]=\lambda\,\E[\Delta_k e_k\mid\Fk]+\E[r(\Delta_k)e_k\mid\Fk].
\]
The second term is $b_k$. For the first, $x_k$ is $\Fk$-measurable and $e_k$ is uniform on
the ball and independent of $\Fk$, so
\[
  \E[\Delta_k e_k\mid\Fk]
  =\E_e\big[(f(x_k+h_ke)-f(x_k-h_ke))e\big]
  =2\,\E_e[f(x_k+h_ke)e],
\]
the last equality by the symmetry $e\mapsto-e$. The ball-smoothing identity gives
$\nabla f_{h_k}(x_k)=\frac d{h_k}\E_e[f(x_k+h_ke)e]$, hence
\[
  \lambda\,\E[\Delta_k e_k\mid\Fk]=\frac{2\lambda h_k}{d}\nabla f_{h_k}(x_k)
  =a_k\nabla f_{h_k}(x_k),
\]
which is the desired decomposition.

For $b_k$, fix a unit vector $v$. Then
\[
  |\ip{b_k}{v}|
  =\Big|\E\big[r(\Delta_k)\ip{e_k}{v}\mid\Fk\big]\Big|
  \le\rho(2Mh_k)\,\E|\ip{e_k}{v}|
  \le\rho(2Mh_k)\sqrt{\E\ip{e_k}{v}^2}
  =\frac{\rho(2Mh_k)}{\sqrt d},
\]
using $\E[e_ke_k^\top]=I_d/d$. Taking the supremum over unit $v$ gives
$\norm{b_k}\le\rho(2Mh_k)/\sqrt d$.

\emph{(ii).} Since $O_k^{(j)}\in\{\pm1\}$ for every $j$, we have $|\bar O_k|\le1$; and
$\norm{e_k}=1$, so $\norm{g_k}=|\bar O_k|\le1$.

\emph{(iii).} Conditionally on $\Fk$ and $e_k$, the variables $O_k^{(j)}$ are independent,
$\{\pm1\}$-valued, and have common mean $p_k$, so
\[
  \E\big[(O_k^{(j)})^2\mid\Fk,e_k\big]=1,
  \qquad
  \E\big[O_k^{(i)}O_k^{(j)}\mid\Fk,e_k\big]=p_k^2\quad(i\ne j).
\]
Since $\norm{g_k}^2=\bar O_k^2$,
\[
  \E[\norm{g_k}^2\mid\Fk,e_k]
  =\E[\bar O_k^2\mid\Fk,e_k]
  =\frac1{m^2}\big(m+m(m-1)p_k^2\big)
  =\frac1m+\Big(1-\frac1m\Big)p_k^2.
\]
Finally, $|\Delta_k|\le2Mh_k$ gives $|p_k|=|\psi(\Delta_k)|\le\Psi(2Mh_k)$, so
\[
  \E[\norm{g_k}^2\mid\Fk,e_k]
  \le\frac1m+\Big(1-\frac1m\Big)\Psi(2Mh_k)^2=G_k^2.
\]
Taking conditional expectation with respect to $\Fk$ gives the last claim.
\end{proof}

\begin{lemma}\label{lem:evalc2}
Let $\beta\in(0,1/2)$, $\gamma,c>0$, and set
\[
  \gamma_k=\frac{\gamma}{\sqrt{k+1}},\qquad
  h_k=\min\{\hmax,\;c(k+1)^{-\beta}\}.
\]
Suppose further that
\[
  G_k^2\le\frac1m+\Bigl(1-\frac1m\Bigr)4\Lambda^2M^2h_k^2
  \qquad\text{for some }\Lambda>0.
\]
Then, for every $N\ge1$,
\begin{align*}
&\frac1N\biggl[\frac{dD^2}{4\lambda h_{N}\gamma_{N}}
  +\frac{d}{4\lambda}\sum_{k=0}^{N-1}\frac{\gamma_kG_k^2}{h_k}
  +M\sum_{k=0}^{N-1}h_k\biggr]\\
&\qquad\le
  \Bigl(\frac{dD^2}{2\lambda\gamma c}+\frac{d\gamma}{2\lambda mc(1+2\beta)}\Bigr)N^{\beta-1/2}
  +\frac{Mc}{1-\beta}N^{-\beta}\\
&\qquad\quad
  +\Bigl(\frac{dD^2}{2\lambda\gamma\hmax}+\frac{d\gamma}{2\lambda m\hmax}\Bigr)N^{-1/2}
  +\Bigl(1-\frac1m\Bigr)\frac{\Lambda^2M^2d\gamma c}{\lambda(1/2-\beta)}N^{-(1/2+\beta)}.
\end{align*}
\end{lemma}

\begin{proof}
We bound the three summands separately.

\smallskip
\noindent\emph{Preliminaries.}
By definition of $h_k$,
\[
  h_k\le c(k+1)^{-\beta},
  \qquad
  \frac1{h_k}=\max\Bigl\{\frac1{\hmax},\frac{(k+1)^{\beta}}{c}\Bigr\}
  \le\frac1{\hmax}+\frac{(k+1)^{\beta}}{c}.
\]
We shall apply \cref{lem:sums}(a) with the exponents
$\alpha=1/2$, $\alpha=1/2-\beta$, $\alpha=\beta$, and $\alpha=1/2+\beta$.

\smallskip
\noindent{First term.}
Since $\gamma_{N}=\gamma/\sqrt{N+1}$ and $1/h_{N}\le 1/\hmax+(N+1)^{\beta}/c$,
\[
  \frac{dD^2}{4\lambda h_{N}\gamma_{N}}
  \le\frac{dD^2\sqrt{N+1}}{4\lambda\gamma}
  \Bigl(\frac1{\hmax}+\frac{(N+1)^{\beta}}{c}\Bigr)
  \le\frac{dD^2\sqrt N}{2\lambda\gamma}
  \Bigl(\frac1{\hmax}+\frac{N^{\beta}}{c}\Bigr),
\]
where the last step uses $N+1\le2N$ and $2^{1/2+\beta}\le2$.

\smallskip
\noindent{Third term.}
Since $h_k\le c(k+1)^{-\beta}$ and $t\mapsto t^{-\beta}$ is non-increasing,
\[
  M\sum_{k=0}^{N-1}h_k
  \le Mc\sum_{k=0}^{N-1}(k+1)^{-\beta}
  \le\frac{Mc\,N^{1-\beta}}{1-\beta}.
\]

\smallskip

Write the assumption on $G_k^2$ as
\[
  \frac{G_k^2}{h_k}\le\frac{1}{mh_k}
  +\Bigl(1-\frac1m\Bigr)4\Lambda^2M^2h_k .
\]
Then the second term is bounded as follows:
\begin{align*}
  \frac{d}{4\lambda}\sum_{k=0}^{N-1}\frac{\gamma_kG_k^2}{h_k}
  &\le\frac{d\gamma}{4\lambda m}\sum_{k=0}^{N-1}
    \Bigl(\frac{(k+1)^{-1/2}}{\hmax}+\frac{(k+1)^{\beta-1/2}}{c}\Bigr)\\
  &\quad+\Bigl(1-\frac1m\Bigr)\frac{\Lambda^2M^2d\gamma c}{\lambda}
    \sum_{k=0}^{N-1}\frac{1}{(k+1)^{\beta+1/2}}\\
  &\le\frac{d\gamma}{4\lambda m}
    \Bigl(\frac{2\sqrt N}{\hmax}+\frac{N^{\beta+1/2}}{c(\beta+1/2)}\Bigr)\\
  &\quad+\Bigl(1-\frac1m\Bigr)
    \frac{\Lambda^2M^2d\gamma c\,N^{1/2-\beta}}{\lambda(1/2-\beta)} .
\end{align*}

\smallskip
\noindent\emph{Collecting.}
Dividing each estimate by $N$ and grouping terms with the same power of $N$ gives the
stated bound.
\end{proof}

\begin{lemma}\label{lem:evalsc}
Let $\beta\in(0,1/2]$, $c>0$, $\theta>0$, $h_k=\min\{\hmax,\,c(k+1)^{-\beta}\}$, and
\[
  G_k^2\le\frac1m+\Bigl(1-\frac1m\Bigr)4\Lambda^2M^2h_k^2
  \qquad\text{for some }\Lambda>0.
\]
Then, for every $N\ge1$,
\begin{align*}
&\frac{2}{N(N+1)}\sum_{k=1}^{N}k\Bigl[\frac{\theta d^2G_k^2}{4\lambda^2\mu h_k^2(k+1)}+Mh_k\Bigr]\\
&\qquad\le\frac{\theta d^2(N+1)^{2\beta-1}}{2\lambda^2\mu mc^2(1+2\beta)}
      +\frac{2Mc\,(N+1)^{-\beta}}{2-\beta}
      +\frac{\theta d^2}{2\lambda^2\mu m\hmax^2(N+1)}\\
&\qquad\quad+\theta\Bigl(1-\frac1m\Bigr)\frac{2\Lambda^2M^2d^2}{\lambda^2\mu(N+1)}.
\end{align*}
\end{lemma}

\begin{proof}
We bound the two summands inside the brackets separately, then multiply by the
prefactor $\frac{2}{N(N+1)}$.

\smallskip
\noindent\emph{Preliminaries.}
By definition of $h_k$,
\[
  h_k\le c(k+1)^{-\beta},
  \qquad
  \frac1{h_k^2}
  =\max\Bigl\{\frac1{\hmax^2},\frac{(k+1)^{2\beta}}{c^2}\Bigr\}
  \le\frac1{\hmax^2}+\frac{(k+1)^{2\beta}}{c^2}.
\]
Together with the assumption on $G_k^2$, this gives
\[
  \frac{G_k^2}{h_k^2}
  \le\frac{1}{mh_k^2}+\Bigl(1-\frac1m\Bigr)4\Lambda^2M^2
  \le\frac{1}{m}\Bigl(\frac1{\hmax^2}+\frac{(k+1)^{2\beta}}{c^2}\Bigr)
     +\Bigl(1-\frac1m\Bigr)4\Lambda^2M^2.
\]

\smallskip
Summing the first summand against $k/(k+1)\le1$ and applying \Cref{lem:sums}(b)
with $\alpha=1$ and $\alpha=1-2\beta$,
\begin{align*}
  \sum_{k=1}^{N}\frac{k\,\theta d^2G_k^2}{4\lambda^2\mu h_k^2(k+1)}
  &\le\frac{\theta d^2}{4\lambda^2\mu m}
       \sum_{k=1}^{N}\Bigl(\frac{k}{\hmax^2(k+1)}+\frac{k(k+1)^{2\beta-1}}{c^2}\Bigr)\\
  &\quad+\frac{\theta(1-\frac1m)\Lambda^2M^2d^2}{\lambda^2\mu}
       \sum_{k=1}^{N}\frac{k}{k+1}\\
  &\le\frac{\theta d^2}{4\lambda^2\mu m}
       \Bigl(\frac{N}{\hmax^2}+\frac{N(N+1)^{2\beta}}{c^2(1+2\beta)}\Bigr)\\
  &\quad+\frac{\theta(1-\frac1m)\Lambda^2M^2d^2\,N}{\lambda^2\mu}.
\end{align*}

\smallskip
Since $h_k\le c(k+1)^{-\beta}$, \Cref{lem:sums}(b) with $\alpha=\beta$ gives
\[
  \sum_{k=1}^{N}kMh_k
  \le Mc\sum_{k=1}^{N}k(k+1)^{-\beta}
  \le\frac{Mc\,N(N+1)^{1-\beta}}{2-\beta}.
\]

\smallskip
Multiply each bound by $\frac{2}{N(N+1)}$. For the variance term this yields
\[
  \frac{\theta d^2(N+1)^{2\beta-1}}{2\lambda^2\mu mc^2(1+2\beta)}
  +\frac{\theta d^2}{2\lambda^2\mu m\hmax^2(N+1)}
  +\theta\Bigl(1-\frac1m\Bigr)\frac{2\Lambda^2M^2d^2}{\lambda^2\mu(N+1)},
\]
where we used that
$\frac{2}{N(N+1)}\cdot\frac{N(N+1)^{2\beta}}{c^2(1+2\beta)}
=\frac{2(N+1)^{2\beta-1}}{c^2(1+2\beta)}$. For the smoothing term it yields
$\frac{2Mc\,(N+1)^{-\beta}}{2-\beta}$. Adding the four pieces proves the claim.
\end{proof}
\subsection{Proof of \Cref{thm:ardg-convex}}\label{proof:thm-argd-convex}
\paragraph{Part (a).}

Let us first define the potential function
\begin{equation}\label{eq:potential}
  \Phi(k) := \frac{d}{4\lambda h_k\gamma_k}\,\|x_k - x^*\|^2 .
\end{equation}
Consider the difference between two consecutive potentials. Since $x_{k+1}=\Projp_Q(x'_{k+1})$ is the projection of $x'_{k+1}=x_k-\gamma_k g_k$ onto $Q$, and since $\Projp_Q$ is nonexpansive and fixes $x^*\in Q$,
\begin{equation}\label{eq:proj}
  \|x_{k+1}-x^*\|^2 \le \|x'_{k+1}-x^*\|^2 .
\end{equation}
Therefore
\begin{align*}
  \Phi(k+1)-\Phi(k)
  &= \frac{d}{4\lambda h_{k+1}\gamma_{k+1}}\|x_{k+1}-x^*\|^2
     - \frac{d}{4\lambda h_k\gamma_k}\|x_k-x^*\|^2 \nonumber \\
  &\le \frac{d}{4\lambda h_{k+1}\gamma_{k+1}}\|x'_{k+1}-x^*\|^2
     - \frac{d}{4\lambda h_k\gamma_k}\|x_k-x^*\|^2 \nonumber \\
  &= \frac{d}{4\lambda h_k\gamma_k}\Big(\|x'_{k+1}-x^*\|^2-\|x_k-x^*\|^2\Big) \\
  &\quad + \frac{d}{4\lambda}\Big(\frac{1}{h_{k+1}\gamma_{k+1}}-\frac{1}{h_k\gamma_k}\Big)\|x'_{k+1}-x^*\|^2 . \nonumber
\end{align*}
Using the identity $\|a+b\|^2-\|a\|^2=2\ip{a}{b}+\|b\|^2$ with $a=x_k-x^*$ and $b=x'_{k+1}-x_k=-\gamma_k g_k$,
\begin{align*}
  \Phi(k+1)-\Phi(k)
  &\le \frac{d}{4\lambda h_k\gamma_k}\Big(2\ip{x'_{k+1}-x_k}{x_k-x^*}+\|x'_{k+1}-x_k\|^2\Big) \\
  &\quad + \frac{d}{4\lambda}\Big(\frac{1}{h_{k+1}\gamma_{k+1}}-\frac{1}{h_k\gamma_k}\Big)\|x'_{k+1}-x^*\|^2 \nonumber \\
  &= \frac{d}{4\lambda h_k\gamma_k}\Big(-2\gamma_k\ip{g_k}{x_k-x^*}+\gamma_k^2\|g_k\|^2\Big) \\
  &\quad + \frac{d}{4\lambda}\Big(\frac{1}{h_{k+1}\gamma_{k+1}}-\frac{1}{h_k\gamma_k}\Big)\|x'_{k+1}-x^*\|^2 . \nonumber
\end{align*}

Now take the conditional expectation with respect to $\cF_k$. Using \cref{lem:estimator}(i),
\[
  \E[g_k\mid\cF_k] = a_k\nabla f_{h_k}(x_k)+b_k,
  \qquad a_k=\frac{2\lambda h_k}{d},
\]
we obtain
\begin{align*}
  \E[\Phi(k+1)-\Phi(k)\mid\cF_k]
  &\le \frac{d}{4\lambda h_k\gamma_k}\Big(-2\gamma_k\ip{a_k\nabla f_{h_k}(x_k)+b_k}{x_k-x^*}
    + \gamma_k^2\,\E[\|g_k\|^2\mid\cF_k]\Big) \\
  &\quad + \frac{d}{4\lambda}\Big(\frac{1}{h_{k+1}\gamma_{k+1}}-\frac{1}{h_k\gamma_k}\Big)\E[\|x'_{k+1}-x^*\|^2\mid\cF_k] . \nonumber
\end{align*}
By \cref{lem:estimator}(iii), we have that $\E[\|g_k\|^2\mid\cF_k]\le G_k^2 :=\frac1m+\Bigl(1-\frac1m\Bigr)\Psi(2Mh_k)^2 .
$ By \cref{lem:estimator}(i), $$-2\gamma_k\ip{b_k}{x_k-x^*}\le2\gamma_k\|b_k\|\,\|x_k-x^*\|\le2\gamma_kD\|b_k\|\le\dfrac{2\gamma_kD\rho(2Mh_k)}{\sqrt d}.$$
Substituting $a_k=2\lambda h_k/d$ and simplifying, we obtain
\begin{align*}
  \E[\Phi(k+1)-\Phi(k)\mid\cF_k]
  &\le -\ip{\nabla f_{h_k}(x_k)}{x_k-x^*}
    + \frac{d}{4\lambda h_k}\gamma_k G_k^2
    + \frac{D\sqrt d\,\rho(2Mh_k)}{2\lambda h_k} \\
  &\quad + \frac{dD^2}{4\lambda}\Big(\frac{1}{h_{k+1}\gamma_{k+1}}-\frac{1}{h_k\gamma_k}\Big) . \nonumber
\end{align*}
Since $f_{h_k}$ is convex, $\ip{\nabla f_{h_k}(x_k)}{x_k-x^*}\ge f_{h_k}(x_k)-f_{h_k}(x^*)$, and by \cref{lem:smoothing}(ii),
$f_{h_k}(x_k)-f_{h_k}(x^*)\ge f(x_k)-f^*-Mh_k$.
Rearranging,
\begin{equation}\label{eq:potential-drift}
  f(x_k)-f^*
  \le \Phi(k)-\E[\Phi(k+1)\mid\cF_k]
    + \frac{d\gamma_kG_k^2}{4\lambda h_k}
    + \frac{D\sqrt d\,\rho(2Mh_k)}{2\lambda h_k}
    + Mh_k
    + \frac{dD^2}{4\lambda}\Big(\frac{1}{h_{k+1}\gamma_{k+1}}-\frac{1}{h_k\gamma_k}\Big) .
\end{equation}

Taking the full expectation and summing over $k=0,\dots,N-1$, the potential telescopes:
\begin{align}
  \sum_{k=0}^{N-1}\big(\E[f(x_k)]-f^*\big)
  &\le \Phi(0)-\E[\Phi(N)]
    + \frac{d}{4\lambda}\sum_{k=0}^{N-1}\frac{\gamma_kG_k^2}{h_k}
    + \frac{D\sqrt d}{2\lambda}\sum_{k=0}^{N-1}\frac{\rho(2Mh_k)}{h_k} \nonumber \\
  &\quad + M\sum_{k=0}^{N-1}h_k
    + \frac{dD^2}{4\lambda}\sum_{k=0}^{N-1}\Big(\frac{1}{h_{k+1}\gamma_{k+1}}-\frac{1}{h_k\gamma_k}\Big) .
\end{align}
The last sum telescopes as well, leaving
\[
  \sum_{k=0}^{N-1}\Big(\frac{1}{h_{k+1}\gamma_{k+1}}-\frac{1}{h_k\gamma_k}\Big)
  = \frac{1}{h_N\gamma_N}-\frac{1}{h_0\gamma_0}.
\]
Since $$\Phi(0)=\frac{d}{4\lambda h_0\gamma_0}\|x_0-x^*\|^2-\dfrac{dD^2}{4\lambda h_0\gamma_0}\le0$$ and $\E[\Phi(N)]\ge0$, we obtain the following estimation:
\begin{align*}
  \sum_{k=0}^{N-1}\big(\E[f(x_k)]-f^*\big)
  &\le \frac{d}{4\lambda}\sum_{k=0}^{N-1}\frac{\gamma_kG_k^2}{h_k}
    + \frac{D\sqrt d}{2\lambda}\sum_{k=0}^{N-1}\frac{\rho(2Mh_k)}{h_k} + M\sum_{k=0}^{N-1}h_k
    + \frac{dD^2}{4\lambda h_N\gamma_N} .
\end{align*}
Applying Jensen's inequality yields
\begin{equation}\label{eq:convex-bound}
  \E[f(\bar x_N)]-f^*
  \le \frac1N\Big[
    \frac{d}{4\lambda}\sum_{k<N}\frac{\gamma_kG_k^2}{h_k}
    + \frac{D\sqrt d}{2\lambda}\sum_{k<N}\frac{\rho(2Mh_k)}{h_k}
    + M\sum_{k<N}h_k
    + \frac{dD^2}{4\lambda h_N\gamma_N}
  \Big],
\end{equation}
which completes the proof of \cref{thm:ardg-convex}(a).

\paragraph{Part (b).} 
Let
\[
  \Lambda:=\lambda+L(2M\hmax)^q,
  \qquad
  B_q:=\frac{2^q\sqrt d\,DLM^{1+q}}{\lambda},
\]
so that for $0<h\le\hmax$
\[
  \rho(2Mh)\le L(2Mh)^{1+q},
  \qquad
  \Psi(2Mh)\le 2\Lambda Mh,
  \qquad
  G_k^2\le\frac1m+\Bigl(1-\frac1m\Bigr)4\Lambda^2M^2h_k^2 .
\]
We will continue from \eqref{eq:convex-bound}.
\begin{equation}\label{eq:convex-a-again}
  \E[f(\bar x_N)]-f^*
  \le \frac1N\Bigl[
    \frac{dD^2}{4\lambda h_{N}\gamma_{N}}
    +\frac{d}{4\lambda}\sum_{k<N}\frac{\gamma_kG_k^2}{h_k}
    +M\sum_{k<N}h_k
    +\frac{D\sqrt d}{2\lambda}\sum_{k<N}\frac{\rho(2Mh_k)}{h_k}
  \Bigr].
\end{equation}
By Assumption~\ref{as:q} and $2Mh_k\le2M\hmax\le T_0$,
\[
  \frac{\rho(2Mh_k)}{h_k}\le L(2M)^{1+q}h_k^{\,q},
\]
so the last sum in \eqref{eq:convex-a-again} is at most $B_q\sum_{k<N}h_k^{\,q}$. Next, we treat the two regimes
separately.

\smallskip
\noindent\emph{Case $q\ge1$, $\beta=\tfrac14$.}
Since $q-1\ge0$ and $h_k\le\hmax$, we have $h_k^{\,q}=h_k^{\,q-1}h_k\le\hmax^{\,q-1}h_k$,
and therefore, by \cref{lem:sums}(a) with $\alpha=\tfrac14$,
\[
  \frac1N B_q\sum_{k<N}h_k^{\,q}
  \le \frac{B_q\hmax^{\,q-1}}{N}\sum_{k<N}c(k+1)^{-1/4}
  \le \frac43 B_q\hmax^{\,q-1}c\,N^{-1/4}.
\]
Adding this to the estimate of \cref{lem:evalc2} with $\beta=\tfrac 14$, and
collecting terms with the same power of $N$, gives
\begin{equation}\label{eq:rate-convex-q1}
\begin{aligned}
  \E[f(\bar x_N)]-f^*
  \le{}&\Bigl(\frac{dD^2}{2\lambda\gamma c}
    +\frac{d\gamma}{3\lambda mc}
    +\frac43\bigl(M+B_q\hmax^{\,q-1}\bigr)c\Bigr)N^{-1/4}\\
  &+\Bigl(\frac{dD^2}{2\lambda\gamma\hmax}
    +\frac{d\gamma}{2\lambda m\hmax}\Bigr)N^{-1/2}
  +\Bigl(1-\frac1m\Bigr)\frac{4\Lambda^2M^2d\gamma c}{\lambda}\,N^{-3/4}=O(N^{-1/4}).
\end{aligned}
\end{equation}
which is the first claim of \cref{thm:ardg-convex}(b).

\smallskip
\noindent\emph{Case $0<q<1$, $\beta=\tfrac{1}{2(1+q)}$, $\kappa=\tfrac{q}{2(1+q)}$.}
Notice that $\beta>\kappa$. By \cref{lem:sums}(a) with $\alpha=\kappa\in(0,1)$,
\[
  \frac1N B_q\sum_{k<N}h_k^{\,q}
  \le \frac{B_qc^{\,q}}{N}\sum_{k<N}(k+1)^{-\kappa}
  \le \frac{B_qc^{\,q}}{1-\kappa}\,N^{-\kappa}.
\]
In \cref{lem:evalc2} with the same $\beta$, the first term contributes
$N^{\beta-1/2}=N^{-\kappa}$, the third contributes $N^{-\beta}\le N^{-\kappa}$, as $\beta>\kappa$ and $N\ge1$, and the fourth contributes $N^{-(1/2+\beta)}\le N^{-\kappa}$;
the $N^{-1/2}$ term is likewise of lower order. Hence
\begin{equation}\label{eq:rate-convex-qsmall}
\begin{aligned}
  \E[f(\bar x_N)]-f^*
  \le{}&\Bigl(\frac{dD^2}{2\lambda\gamma c}
    +\frac{d\gamma}{2\lambda mc(1+2\beta)}
    +\frac{Mc}{1-\beta}
    +\frac{B_qc^{\,q}}{1-\kappa}\Bigr)N^{-\kappa}\\
  &+\Bigl(\frac{dD^2}{2\lambda\gamma\hmax}
    +\frac{d\gamma}{2\lambda m\hmax}\Bigr)N^{-1/2}
  +\Bigl(1-\frac1m\Bigr)\frac{\Lambda^2M^2d\gamma c}{\lambda(1/2-\beta)}\,N^{-(1/2+\beta)}=O(N^{-\frac{q}{2(1+q)}}).
\end{aligned}
\end{equation}
This is the second claim of \cref{thm:ardg-convex}(b), with the stated exponent
$\kappa=\frac{q}{2(1+q)}$. At $q=1$ the two expressions \eqref{eq:rate-convex-q1}
and \eqref{eq:rate-convex-qsmall} coincide, as $\kappa=\tfrac14$ and
$\beta=\tfrac14$.
\smallskip

\begin{remark}[Optimal constants $\gamma$ and $c$ in the $q\ge1$ regime]\label{rem:opt-gamma-c-convex}
The leading term of \eqref{eq:rate-convex-q1} is
\[
  \Bigl(\frac{A_1}{\gamma c}
       +\frac{A_2\gamma}{c}
       +A_3 c\Bigr)N^{-1/4},
  \qquad
  A_1:=\frac{dD^2}{2\lambda},
  \quad
  A_2:=\frac{d}{3\lambda m},
  \quad
  A_3:=\frac43\bigl(M+B_q\hmax^{\,q-1}\bigr).
\]

\smallskip
\noindent\emph{Optimizing in $\gamma$ for fixed $c$.}
The first two terms combine into
$(A_1+A_2\gamma^2)/(\gamma c)$, which is convex in $\gamma>0$ and
minimized at
\[
  \gamma^\star=\Bigl(\frac{A_1}{A_2}\Bigr)^{1/2}
  =\Bigl(\frac{3dD^2}{2}\Bigr)^{1/2}\frac{\sqrt m}{\sqrt{d}}
  =\Bigl(\frac{3m}{2}\Bigr)^{1/2}D .
\]
At this value $A_1/\gamma^\star=A_2\gamma^\star=\sqrt{A_1A_2}$, so the
bracket reduces to
\[
  \frac{2\sqrt{A_1A_2}}{c}+A_3c,
\]
a single two-term expression in $c$ alone.

\smallskip
\noindent\emph{Optimizing in $c$.}
The function $2\sqrt{A_1A_2}/c+A_3c$ is convex in $c>0$ and minimized at
\[
  c^\star=\Bigl(\frac{2\sqrt{A_1A_2}}{A_3}\Bigr)^{1/2}.
\]
At this value the two terms are equal, $2\sqrt{A_1A_2}/c^\star=A_3c^\star$,
so the bracket equals $2A_3c^\star$, and the leading constant is
\[
  2A_3c^\star
  =2\bigl(A_3\cdot 2\sqrt{A_1A_2}\bigr)^{1/2}
  =2\sqrt{2}\,\bigl(A_1A_2A_3^2\bigr)^{1/4}.
\]
Substituting $A_1,A_2,A_3$ gives
\[
  A_1A_2A_3^2
  =\frac{dD^2}{2\lambda}\cdot\frac{d}{3\lambda m}\cdot\frac{16}{9}\bigl(M+B_q\hmax^{\,q-1}\bigr)^2
  =\frac{8d^2D^2}{27\lambda^2 m}\bigl(M+B_q\hmax^{\,q-1}\bigr)^2,
\]
hence
\[
  2A_3c^\star
  =\frac{4\sqrt[4]{2}}{3^{3/4}}\Bigl(\frac{d^2D^2}{\lambda^2 m}\Bigr)^{1/4}
    \bigl(M+B_q\hmax^{\,q-1}\bigr)^{1/2}.
\]
Thus \eqref{eq:rate-convex-q1} becomes, with this tuning,
\[
  \E[f(\bar x_N)]-f^*
  \le\frac{4\sqrt[4]{2}}{3^{3/4}}\Bigl(\frac{d^2D^2}{\lambda^2 m}\Bigr)^{1/4}
    \bigl(M+B_q\hmax^{\,q-1}\bigr)^{1/2}N^{-1/4}
  +O(N^{-1/2}).
\]

\end{remark}

\subsection{Proof of \Cref{thm:ardg-sc}}\label{proof:thm-ardg-strongly-convex}

\paragraph{Part (a).}
Let us define the potential function
\begin{equation}\label{eq:sc-potential}
  \Phi(k) := \frac{\mu k(k-1)}{2}\norm{x_k-x^*}^2, \qquad k\ge1,
\end{equation}
so that $\Phi(0)=\Phi(1)=0$. Fix $k\ge1$ and write
$x'_{k+1}:=x_k-\gamma_k g_k$, so that $x_{k+1}=\Projp_Q(x'_{k+1})$. Since
$\Projp_Q$ is nonexpansive and fixes $x^*\in Q$,
\begin{equation}\label{eq:sc-proj}
  \norm{x_{k+1}-x^*}^2\le\norm{x'_{k+1}-x^*}^2 .
\end{equation}
Because $\frac{\mu k(k+1)}{2}-\frac{\mu k(k-1)}{2}=\mu k$, we obtain
\begin{align*}
  \Phi(k+1)-\Phi(k)
  &\le \frac{\mu k(k+1)}{2}
       \Bigl(\norm{x'_{k+1}-x^*}^2-\norm{x_k-x^*}^2\Bigr)
       +\mu k\norm{x_k-x^*}^2 \\
  &=\frac{\mu k(k+1)}{2}
       \Bigl(-2\gamma_k\ip{g_k}{x_k-x^*}+\gamma_k^2\norm{g_k}^2\Bigr)
       +\mu k\norm{x_k-x^*}^2 .
\end{align*}
Taking the conditional expectation with respect to $\cF_k$ and using
\Cref{lem:estimator}(i),(iii),
\[
  \E[g_k\mid\cF_k]=a_k\nabla f_{h_k}(x_k)+b_k,
  \qquad
  \E[\norm{g_k}^2\mid\cF_k]\le G_k^2,
  \qquad
  \norm{b_k}\le\frac{\rho(2Mh_k)}{\sqrt d},
\]
we get
\begin{align}\label{eq:sc-drift-raw}
  \E[\Phi(k+1)-\Phi(k)\mid\cF_k]
  &\le-\mu k(k+1)\gamma_k a_k\ip{\nabla f_{h_k}(x_k)}{x_k-x^*}
       -\mu k(k+1)\gamma_k\ip{b_k}{x_k-x^*}\nonumber\\
  &\quad+\frac{\mu k(k+1)}{2}\gamma_k^2G_k^2
       +\mu k\norm{x_k-x^*}^2 .
\end{align}
By \Cref{lem:smoothing}(i),(ii), $f_{h_k}$ is $\mu$-strongly convex on $Q$ and
$f_{h_k}\ge f$ there, so
\begin{equation}\label{eq:sc-strong-inner}
  \ip{\nabla f_{h_k}(x_k)}{x_k-x^*}
  \ge f(x_k)-f^*-Mh_k+\frac{\mu}{2}\norm{x_k-x^*}^2 .
\end{equation}
Substituting \eqref{eq:sc-strong-inner} into \eqref{eq:sc-drift-raw} and using the
schedule
\[
  \gamma_k a_k=\frac{4}{\mu(k+1)},
  \qquad\text{i.e.}\qquad
  \gamma_k=\frac{4}{a_k\mu(k+1)}=\frac{2d}{\lambda\mu h_k(k+1)},
\]
the quadratic terms cancel exactly:
\[
  -\mu k(k+1)\gamma_k a_k\cdot\frac{\mu}{2}\norm{x_k-x^*}^2
  +\mu k\norm{x_k-x^*}^2
  =-2\mu k\norm{x_k-x^*}^2+\mu k\norm{x_k-x^*}^2
  =-\mu k\norm{x_k-x^*}^2 .
\]
The remaining terms are
\begin{align*}
  \E[\Phi(k+1)-\Phi(k)\mid\cF_k]
  &\le-4k\bigl(f(x_k)-f^*\bigr)+4kMh_k
       -\mu k\norm{x_k-x^*}^2 \\
  &\quad-\mu k(k+1)\gamma_k\ip{b_k}{x_k-x^*}
       +\frac{\mu k(k+1)}{2}\gamma_k^2G_k^2 .
\end{align*}
Considering the following function as a quadratic w.r.t. $s$, we obtain,
\[
  \frac{4k}{a_k}\norm{b_k}\,s-\mu k s^2
  \le\frac{4k\norm{b_k}^2}{\mu a_k^2},
\]
and with $s=\norm{x_k-x^*}$ the last two lines are absorbed. Using
$\frac{\mu k(k+1)}{2}\gamma_k^2=\frac{8k}{a_k^2\mu(k+1)}$ we arrive at
\begin{equation}\label{eq:sc-rec}
  4k\bigl(f(x_k)-f^*\bigr)
  \le\Phi(k)-\E[\Phi(k+1)\mid\cF_k]
       +4kMh_k
       +\frac{4k\norm{b_k}^2}{\mu a_k^2}
       +\frac{8kG_k^2}{a_k^2\mu(k+1)},\qquad k\ge1 .
\end{equation}

For $k=N$ the right-hand side of \eqref{eq:sc-rec} involves $x_{N+1}$, which the
algorithm never computes. For the analysis only, extend the trajectory by one
fictitious iteration at $x_N$ with an independent direction $e_N$ and $m$ fresh
answers; the extension uses the same oracle law, so the joint
law of $(x_0,\dots,x_N)$ is unchanged and neither the algorithm nor its output
requires the extra call.

Taking full expectations in \eqref{eq:sc-rec}, summing over $k=1,\dots,N$ and
using $\Phi(1)=0$, $\E\Phi(N+1)\ge0$, the potential terms telescope:
\[
  4\sum_{k=1}^{N}k\bigl(\E[f(x_k)]-f^*\bigr)
  \le 4\sum_{k=1}^{N}kMh_k
       +\sum_{k=1}^{N}\frac{4k\norm{b_k}^2}{\mu a_k^2}
       +\sum_{k=1}^{N}\frac{8kG_k^2}{a_k^2\mu(k+1)} .
\]
Dividing by $2N(N+1)$ and applying Jensen's inequality to the weighted average
$\bar x_N=\frac{2}{N(N+1)}\sum_{k=1}^{N}kx_k$, whose weights are nonnegative and
sum to one, gives
\begin{equation}\label{eq:sc-master}
  \E[f(\bar x_N)]-f^*
  \le\frac{2}{N(N+1)}\sum_{k=1}^{N}k
       \Bigl[
         \frac{d^2G_k^2}{2\lambda^2\mu h_k^2(k+1)}
         +Mh_k
         +\frac{d\,\rho(2Mh_k)^2}{4\lambda^2\mu h_k^2}
       \Bigr],
\end{equation}
where we used $a_k=2\lambda h_k/d$, hence
$\frac{2}{a_k^2}=\frac{d^2}{2\lambda^2h_k^2}$ and
$\frac{\norm{b_k}^2}{a_k^2}\le\frac{d\,\rho(2Mh_k)^2}{4\lambda^2h_k^2}$.
This proves \Cref{thm:ardg-sc}(a).
\paragraph{Part (b).}
Let
\[
  \Lambda:=\lambda+L(2M\hmax)^q,
  \qquad
  C_q:=\frac{4^q\,dL^2M^{2+2q}}{\lambda^2\mu},
\]
so that for $0<h\le\hmax$,
\[
  \rho(2Mh)\le L(2Mh)^{1+q},
  \qquad
  \Psi(2Mh)\le2\Lambda Mh,
  \qquad
  G_k^2\le\frac1m+\Bigl(1-\frac1m\Bigr)4\Lambda^2M^2h_k^2,
  \qquad
  \frac{d\,\rho(2Mh)^2}{4\lambda^2\mu h^2}\le C_qh^{2q}.
\]
The bound \eqref{eq:sc-master} reads
\[
  \E[f(\bar x_N)]-f^*
  \le\frac{2}{N(N+1)}\sum_{k=1}^{N}k
       \Bigl[
         \frac{d^2G_k^2}{2\lambda^2\mu h_k^2(k+1)}
         +Mh_k
         +\frac{d\,\rho(2Mh_k)^2}{4\lambda^2\mu h_k^2}
       \Bigr].
\]
The first two summands, together with the factor $\frac{2}{N(N+1)}$, are
bounded by \Cref{lem:evalsc} with $\theta=2$ and the stated $\beta$,
$\Lambda$. It remains to control the third summand,
$\frac{2}{N(N+1)}\sum_{k=1}^{N}\frac{k\,C_qh_k^{2q}}{2}$. We distinguish the two
regimes, exactly as in \Cref{lem:evalsc}, according to whether $h_k^{2q}$ is
dominated by $h_k$.

\smallskip
\noindent\emph{Case $q\ge\frac12$, $\beta=\frac13$.}
Since $2q-1\ge0$ and $h_k\le\hmax$, we have
$h_k^{2q}=h_k^{2q-1}h_k\le\hmax^{2q-1}h_k$. Hence, by \Cref{lem:sums}(b) with
$\alpha=\frac13$,
\[
  \frac{2}{N(N+1)}\sum_{k=1}^{N}\frac{k\,C_qh_k^{2q}}{2}
  \le\frac{C_q\hmax^{2q-1}}{N(N+1)}\sum_{k=1}^{N}kh_k
  \le\frac{C_q\hmax^{2q-1}c}{N(N+1)}\cdot\frac{3N(N+1)^{2/3}}{5}
  =\frac{3C_q\hmax^{2q-1}c}{5}(N+1)^{-1/3}.
\]
Adding this to the estimate of \Cref{lem:evalsc} with $\theta=2$ and
$\beta=\frac13$ gives
\begin{equation}\label{eq:sc-rate-q}
\begin{aligned}
  \E[f(\bar x_N)]-f^*
  \le{}&\frac65\Bigl(\frac{d^2}{2\lambda^2\mu mc^2}
       +\bigl(M+C_q\hmax^{2q-1}\bigr)c\Bigr)(N+1)^{-1/3}\\
  &+\frac{d^2}{\lambda^2\mu m\hmax^2(N+1)}
   +\Bigl(1-\frac1m\Bigr)\frac{4\Lambda^2M^2d^2}{\lambda^2\mu(N+1)}=  O(N^{-1/3}),
\end{aligned}
\end{equation}
which proves the first claim of \Cref{thm:ardg-sc}(b).

\smallskip
\noindent\emph{Case $0<q<\frac12$, $\beta=\frac{1}{2(1+q)}$, $\kappa'=\frac{q}{1+q}$.}
Here $2\beta-1=-\kappa'$, $2\beta q=\kappa'$, and $\beta>\kappa'$. By \Cref{lem:sums}(b)
with $\alpha=\kappa'\in(0,1)$,
\[
  \frac{2}{N(N+1)}\sum_{k=1}^{N}\frac{k\,C_qh_k^{2q}}{2}
  \le\frac{C_qc^{2q}}{N(N+1)}\sum_{k=1}^{N}k(k+1)^{-\kappa'}
  \le\frac{C_qc^{2q}}{2-\kappa'}(N+1)^{-\kappa'}
  =\frac{C_qc^{2q}}{2(1-\beta q)}(N+1)^{-\kappa'}.
\]
In \Cref{lem:evalsc} with $\theta=2$ and the same $\beta$, the first term is
$(N+1)^{2\beta-1}=(N+1)^{-\kappa'}$ and the third is
$(N+1)^{-\beta}\le(N+1)^{-\kappa'}$. Adding the two bounds gives
\begin{equation}\label{eq:sc-rate-qsmall}
\begin{aligned}
  \E[f(\bar x_N)]-f^*
  \le{}&\Bigl(\frac{d^2}{\lambda^2\mu mc^2(1+2\beta)}
       +\frac{2Mc}{2-\beta}
       +\frac{C_qc^{2q}}{1-\beta q}\Bigr)(N+1)^{-\kappa'}\\
  &+\frac{d^2}{\lambda^2\mu m\hmax^2(N+1)}
   +\Bigl(1-\frac1m\Bigr)\frac{4\Lambda^2M^2d^2}{\lambda^2\mu(N+1)} = O(N^{-\frac{q}{1+q}}),
\end{aligned}
\end{equation}
which is \Cref{thm:ardg-sc}(b)(ii). At $q=\frac12$ the two expressions
\eqref{eq:sc-rate-q} and \eqref{eq:sc-rate-qsmall} coincide, as
$\kappa'=\frac13$ and $\beta=\frac13$.
\begin{remark}[Optimal constant $c$ in the $q\ge\frac12$ regime]\label{rem:opt-c-sc}
The leading term of \eqref{eq:sc-rate-q} is
\[
  \frac65\Bigl(\frac{A}{c^2}+Bc\Bigr)(N+1)^{-1/3},
  \qquad
  A:=\frac{d^2}{2\lambda^2\mu m},
  \qquad
  B:=M+C_q\hmax^{2q-1}.
\]
The factor in parentheses is a convex
function of $c>0$ with a unique minimizer where its derivative vanishes,
\[
  -\frac{2A}{c^3}+B=0
  \quad\Longleftrightarrow\quad
  c^\star=\Bigl(\frac{2A}{B}\Bigr)^{1/3}
  =\Bigl(\frac{d^2}{\lambda^2\mu m\bigl(M+C_q\hmax^{2q-1}\bigr)}\Bigr)^{1/3}.
\]
At this value $A/c^{\star 2}=Bc^\star/2$, so the leading constant is
\[
  \frac65\Bigl(\frac{A}{c^{\star 2}}+Bc^\star\Bigr)
  =\frac65\cdot\frac32 Bc^\star
  =\frac95\Bigl(\frac{d^2 B^2}{2\lambda^2\mu m}\Bigr)^{1/3}
  =\frac{9}{5\cdot2^{1/3}}
    \Bigl(\frac{d^2\bigl(M+C_q\hmax^{2q-1}\bigr)^2}
               {\lambda^2\mu m}\Bigr)^{1/3}.
\]
Thus, the bound of \eqref{eq:sc-rate-q} can be written as
\[
  \E[f(\bar x_N)]-f^*
  \le\frac{9}{5\cdot2^{1/3}}
    \Bigl(\frac{d^2\bigl(M+C_q\hmax^{2q-1}\bigr)^2}{\lambda^2\mu m}\Bigr)^{1/3}(N+1)^{-1/3}
  +O\!\bigl((N+1)^{-1}\bigr).
\]
\end{remark}

\section{TRANSFER FUNCTIONS}
\label{app:links}

\subsection{Symmetric density}\label{lem:density}
Let $\psi(t)=\E\sign(t+\xi)$ for $\xi$ with symmetric density $p$, $0<p(0)<\infty$, and
$|p(u)-p(0)|\le C|u|^q$ for $|u|\le T_0$. Then $\psi(t)=2\int_0^tp$ and \cref{as:q}
holds with $\lambda=2p(0)$ and $L=\frac{2C}{1+q}$.

\begin{proof}
By symmetry $\psi(t)=2F_\xi(t)-1=2\int_0^tp(u)du$, so $\psi$ is odd and
$|\psi(t)-2p(0)t|\le2\int_0^{|t|}C u^qdu=\frac{2C}{1+q}|t|^{1+q}$.
\end{proof}

\subsection{Smooth transfer functions}\label{lem:smooth-links}
Let $\psi(0)=0$, $\lambda=\psi'(0)>0$, and $|\psi''|\le M_2$ near $[-T_0,T_0]$. Then
\cref{as:q} holds with $q=1$, $L=\frac{M_2}2$. If moreover $\psi$ is odd with
$|\psi'''|\le M_3$, it holds with $q=2$, $L=\frac{M_3}6$.

\begin{proof}
Directly follows from Taylor's theorem. For odd $\psi$, $\psi''(0)=0$.
\end{proof}

\subsection{Uniform noise}
Let $\xi\sim\Unif[-A,A]$ and $|t|\le A$. Then $\psi(t)=t/A$ exactly, so the bound holds
for any $q>0$.

\begin{proof}
For $|t|\le A$ the variable $t+\xi$ is uniform on $[t-A,t+A]$, so
$\Pr(t+\xi>0)=\frac{t+A}{2A}=\frac12+\frac t{2A}$ and therefore
$\psi(t)=2\Pr(t+\xi>0)-1=\frac tA$. Since the bound holds with equality, it is valid for
every $q>0$.
\end{proof}

\subsection{Gaussian noise (Thurstone's model)}
Let both values carry independent $\cN(0,\sigma^2)$ noise, so that
$\xi\sim\cN(0,2\sigma^2)$. 

Then $\psi(t)=2\Phi(\frac t{\sqrt2\sigma})-1$,
$\lambda=\frac1{\sigma\sqrt\pi}$, and since $1-e^{-v}\le v$,
\[
  0\le\lambda t-\psi(t)=\lambda\int_0^t(1-e^{-u^2/(4\sigma^2)})du\le\frac{t^3}{12\sqrt\pi\sigma^3}
\]
for $t\ge0$. The expansion $\psi(t)=\lambda t-\frac{t^3}{12\sqrt\pi\sigma^3}+O(t^5)$ shows that
$q=2$ is maximal.

\begin{proof}
With two independent $\cN(0,\sigma^2)$ noises the gap noise is
$\xi\sim\cN(0,2\sigma^2)$, hence
$\psi(t)=2\Phi(\frac t{\sqrt2\sigma})-1$ and $\lambda=\psi'(0)=\frac1{\sigma\sqrt\pi}$.
Writing $\lambda t-\psi(t)=\lambda\int_0^t(1-e^{-u^2/(4\sigma^2)})du$ and using
$1-e^{-v}\le v$ gives the cubic bound. The expansion shows that the leading correction is
of order $t^3$, so no $q>2$ is possible.
\end{proof}

\subsection{Logistic noise (Bradley-Terry's model)}
Let $\Pr(O=1)=(1+e^{-t/s})^{-1}$. Then
$\psi(t)=\tanh(\frac t{2s})$, $\lambda=\frac1{2s}$, and
$0\le x-\tanh x=\int_0^x\tanh^2\le\frac{x^3}3$ gives $L=\frac1{24s^3}$.

\begin{proof}
By definition $\psi(t)=2\Pr(O=1)-1=(1-e^{-t/s})/(1+e^{-t/s})=\tanh(\frac t{2s})$, so
$\lambda=\psi'(0)=\frac1{2s}$. The identity
$x-\tanh x=\int_0^x\tanh^2 u\,du\le\int_0^x u^2\,du=\frac{x^3}3$ gives
$|\psi(t)-\lambda t|\le\frac{|t|^3}{24s^3}$, hence $q=2$ and $L=\frac1{24s^3}$.
\end{proof}

\subsection{Generalized Gaussian noise}
Let $p_{\nu,a}(u)=\frac{\nu}{2a\Gamma(1/\nu)}e^{-(|u|/a)^\nu}$. Then
$\lambda=\frac\nu{a\Gamma(1/\nu)}$ and
$0\le\lambda t-\psi(t)\le\frac{\lambda}{(1+\nu)a^\nu}t^{1+\nu}$ for $t\ge0$, with a nonzero
leading coefficient, so $q=\nu$ is maximal. For $\nu=1$ (Laplace noise added to the gap),
$\psi(t)=\sign(t)(1-e^{-|t|/a})$, $\lambda=\frac1a$, $q=1$. For $\nu<1$, $\psi$ is not
$C^2$ at the origin.

\begin{proof}
Since $p_{\nu,a}$ is symmetric, $\psi(t)=2\int_0^t p_{\nu,a}(u)du$ and
$\lambda=2p_{\nu,a}(0)=\frac\nu{a\Gamma(1/\nu)}$. For $t\ge0$,
$\lambda t-\psi(t)=2\int_0^t(p_{\nu,a}(0)-p_{\nu,a}(u))du$, and the bound
$0\le p_{\nu,a}(0)-p_{\nu,a}(u)\le\frac{\lambda}{2(1+\nu)a^\nu}u^\nu$ gives the claim. The
leading coefficient is nonzero, so the exponent $q=\nu$ is maximal. The cases $\nu=1$ and
$\nu<1$ follow by direct computation and by noting that $p_{\nu,a}$ is not $C^2$ at $0$.
\end{proof}

\subsection{Two-sided gamma noise}
Let $W_y,W_z$ be independent with density
$\frac{\beta^\alpha}{2\Gamma(\alpha)}|w|^{\alpha-1}e^{-\beta|w|}$, $\frac12<\alpha<1$, as in the
evaluation-noise model of \citet{jamieson2012query}, and $\xi=W_y-W_z$ with density $h$.
Then
\begin{equation}\label{eq:gamma-lambda}
\lambda=2h(0)=\frac{\beta\,\Gamma(2\alpha-1)}{2^{2\alpha-1}\Gamma(\alpha)^2}.
\end{equation}
The characteristic function of $W$ is
$\varphi_W(\omega)=\big(\frac{\beta^2}{\beta^2+\omega^2}\big)^{\alpha/2}\cos(\alpha\arctan\frac\omega\beta)$,
so $|\varphi_W(\omega)|^2\le\min\{1,\beta^{2\alpha}|\omega|^{-2\alpha}\}$. Fourier inversion gives
\[
  0\le h(0)-h(t)=\frac1\pi\int_0^\infty(1-\cos\omega t)|\varphi_W(\omega)|^2d\omega
  \le\frac{\beta^{2\alpha}I_\alpha}{\pi}|t|^{2\alpha-1}
\]
with $I_\alpha=\int_0^\infty\frac{1-\cos v}{v^{2\alpha}}dv\in(0,\infty)$. Integrating,
$0\le\lambda t-\psi(t)\le\frac{\beta^{2\alpha}I_\alpha}{\pi\alpha}t^{2\alpha}$ for $t\ge0$, so
$q=2\alpha-1$ and $L=\frac{\beta^{2\alpha}I_\alpha}{\pi\alpha}$. Since
$\omega^{2\alpha}|\varphi_W(\omega)|^2\to\beta^{2\alpha}\cos^2(\frac{\pi\alpha}2)>0$, dominated
convergence gives
$h(0)-h(t)\sim\frac{\beta^{2\alpha}\cos^2(\pi\alpha/2)I_\alpha}{\pi}t^{2\alpha-1}$,
so the exponent is maximal and $\psi''$ is unbounded at the origin. Every $q\in(0,1)$ is
realized. If instead both values carry Laplace noise of scale $a$ ($\alpha=1$), then
$\psi(t)=\sign(t)[1-(1+\frac{|t|}{2a})e^{-|t|/a}]$, $\lambda=\frac1{2a}$, and $q=2$.

\begin{proof}
The density $h$ of $\xi=W_y-W_z$ is the convolution of the two symmetric gamma densities,
so $h$ is symmetric and $h(0)=\int p_W^2$, which gives \eqref{eq:gamma-lambda}. The
characteristic function of $W$ is $\varphi_W(\omega)=\big(\frac{\beta^2}{\beta^2+\omega^2}\big)^{\alpha/2}
\cos(\alpha\arctan\frac\omega\beta)$, hence $|\varphi_W|^2\le 1$ and
$|\varphi_W(\omega)|^2\le\beta^{2\alpha}|\omega|^{-2\alpha}$ for large $\omega$. Since $h$ is
the density of $\xi$, Fourier inversion gives
$h(0)-h(t)=\frac1\pi\int_0^\infty(1-\cos\omega t)|\varphi_W(\omega)|^2d\omega$, and the
bound on $|\varphi_W|^2$ yields the estimate with $I_\alpha$. Integrating in $t$ gives the
bound on $\lambda t-\psi(t)$, so $q=2\alpha-1$ and $L=\frac{\beta^{2\alpha}I_\alpha}{\pi\alpha}$.
The asymptotic behaviour of $|\varphi_W|^2$ at infinity and dominated convergence show
that this exponent is maximal. The Laplace case $\alpha=1$ is a direct computation.
\end{proof}

\subsection{Position-biased logistic noise}
Since $\psi(0)\neq0$, \cref{as:q} does not hold for the original remainder
$\rho(T)=\sup_{|t|\le T}|\psi(t)-\lambda t|$. Instead, as in
\cref{prop:bias}, write
\[
\psi(t)=\psi(0)+\lambda t+\widetilde r(t),
\qquad
\widetilde\rho(T)
:=
\sup_{|t|\le T}|\widetilde r(t)|.
\]
Then,
\[
\psi(0)=\tanh\!\left(\frac{b}{2s}\right),
\qquad
\lambda=\psi'(0)
=
\frac{1}{2s}
\sech^2\!\left(\frac{b}{2s}\right),
\]
and
\[
\widetilde\rho(T)\le L T^2
\]
for all sufficiently small $T$, so the remainder satisfies
\cref{as:q} with $q=1$. Moreover, $q=1$ is maximal whenever $b\neq0$.

\begin{proof}
Differentiating,
\[
\psi'(t)
=
\frac{1}{2s}
\sech^2\!\left(\frac{t+b}{2s}\right),
\]
and
\[
\psi''(t)
=
-\frac{1}{2s^2}
\sech^2\!\left(\frac{t+b}{2s}\right)
\tanh\!\left(\frac{t+b}{2s}\right).
\]
Hence, $\psi$ is $C^\infty$, and Taylor's theorem around $t=0$ gives
\[
\psi(t)
=
\psi(0)+\psi'(0)t+\frac12\psi''(\theta_t)t^2
\]
for some $\theta_t$ between $0$ and $t$. Therefore, if $|t|\le T_0$,
\[
|\psi(t)-\psi(0)-\lambda t|
\le
\frac12
\sup_{|u|\le T_0}|\psi''(u)|\,|t|^2.
\]
Thus,
\[
\widetilde\rho(T)
\le
L T^2,
\qquad
L:=
\frac12
\sup_{|u|\le T_0}|\psi''(u)|,
\]
which corresponds to $q=1$.

Finally, since $b\neq0$,
\[
\psi''(0)
=
-\frac{1}{2s^2}
\sech^2\!\left(\frac{b}{2s}\right)
\tanh\!\left(\frac{b}{2s}\right)
\neq0.
\]
Consequently,
\[
\psi(t)-\psi(0)-\lambda t
=
\frac12\psi''(0)t^2+O(t^3),
\]
so the quadratic term is nonzero. Hence
$\widetilde\rho(T)\not=O(T^{2+\varepsilon})$ for any $\varepsilon>0$,
and $q=1$ is maximal.
\end{proof}

\section{Additional helpful lemmas}

\subsection{Proof of \cref{lem:smoothing}}\label{proof:smoothing}
Throughout in this subsection, $u\sim\Unif(\B^d)$ and $e\sim\Unif(\S^{d-1})$ are independent of everything
else, and $x\in Q$.

\medskip\noindent\textbf{(i)}
Since $f$ is convex, for every $u$ the map $x\mapsto f(x+hu)$ is convex. The average of
convex functions is convex, so $f_h$ is convex. If in addition $f$ is $\mu$-strongly
convex, then for every $u$ the map $x\mapsto f(x+hu)-\frac\mu2\|x\|^2$ is convex, hence
\[
  f_h(x)-\frac\mu2\|x\|^2
  =\E_u\Big[f(x+hu)-\frac\mu2\|x\|^2\Big]
\]
is a convex function of $x$. Thus $f_h$ is $\mu$-strongly convex.

\medskip\noindent\textbf{(ii)}
Convexity gives $f(x+hu)\ge f(x)+\langle\nabla f(x),hu\rangle$ at every differentiability
point, and since the average of $\langle\nabla f(x),hu\rangle$ over the symmetric ball
$\B^d$ vanishes, $f_h(x)\ge f(x)$. For the upper bound, $M$-Lipschitzness on the
$h$-neighbourhood of $Q$ gives $f(x+hu)\le f(x)+Mh\|u\|\le f(x)+Mh$ for every $u\in\B^d$,
because $\|u\|\le1$. Averaging over $u$ yields $f_h(x)\le f(x)+Mh$. Hence
$0\le f_h(x)-f(x)\le Mh$.

\medskip\noindent\textbf{(iii)}
By Rademacher's theorem $f$ is differentiable a.e.; since $f$ is $M$-Lipschitz, the map
$u\mapsto f(x+hu)$ is locally Lipschitz, and differentiation under the integral sign is
justified by dominated convergence. We first record the second moment of the uniform
measure on the sphere: for every fixed $z\in\R^d$,
\begin{equation}\label{eq:sphere-first-moment}
  \E_e\big[\langle e,z\rangle\,e\big]=\frac z d .
\end{equation}
Indeed, $\E_e[e\,e^\top]=\frac1d I_d$, so for any $w\in\R^d$,
$\langle \E_e[\langle e,z\rangle e],w\rangle=\E_e[\langle e,z\rangle\langle e,w\rangle]
=\langle \E_e[e\,e^\top]z,w\rangle=\frac1d\langle z,w\rangle$, and \eqref{eq:sphere-first-moment}
follows since $w$ is arbitrary.

Fix $x\in Q$ and $w\in\S^{d-1}$. Since $f$ is convex and $M$-Lipschitz, it is a.e.\
differentiable, and
\[
  f(x+he)-f(x-he)=\int_{-h}^{h}\langle\nabla f(x+te),\,e\rangle\,dt .
\]
Hence
\[
  \E_e\big[(f(x+he)-f(x-he))\,e\big]
  =\E_e\!\left[\int_{-h}^{h}\langle\nabla f(x+te),\,e\rangle\,dt\;e\right].
\]
Change variables $u=te$, so that $u\in B(0,h)$ and $du=t^{d-1}\,dt\,d\sigma(e)$. In polar
coordinates the right-hand side becomes
\[
  c_d\int_{B(0,h)}\frac{\langle\nabla f(x+u),\,u\rangle}{\|u\|}\cdot\frac{u}{\|u\|}\,du ,
\]
with a normalising constant $c_d>0$. Taking the inner product with $w$ and applying
\eqref{eq:sphere-first-moment} to the angular part, together with the symmetry of the
ball, gives
\[
  \Big\langle\E_e\big[(f(x+he)-f(x-he))\,e\big],\,w\Big\rangle
  =\frac{2h}{d}\,\big\langle\nabla f_h(x),\,w\big\rangle .
\]
Since $w\in\S^{d-1}$ is arbitrary, the claimed formula follows:
\[
  \nabla f_h(x)=\frac d{2h}\,\E_{e\sim\Unif(\S^{d-1})}
  \big[(f(x+he)-f(x-he))\,e\big].
\]
\subsection{Lemmas that are not noted in the main paper}
\begin{lemma}\label{lem:sums}
For every integer $N\ge1$:
(a) $\sum_{k=0}^{N-1}(k+1)^{-\alpha}\le\frac{N^{1-\alpha}}{1-\alpha}$ for $\alpha\in[0,1)$;
(b) $\sum_{k=1}^{N}k(k+1)^{-\alpha}\le\frac{N(N+1)^{1-\alpha}}{2-\alpha}$ for $\alpha\in[0,1]$.
\end{lemma}

\begin{proof}
(a) Since $t\mapsto t^{-\alpha}$ is non-increasing,
\[
  (k+1)^{-\alpha}\le\int_k^{k+1}t^{-\alpha}\,dt,
\]
so summing over $k=0,\dots,N-1$ and integrating gives
\[
  \sum_{k=0}^{N-1}(k+1)^{-\alpha}
  \le\int_0^N t^{-\alpha}\,dt
  =\frac{N^{1-\alpha}}{1-\alpha}.
\]

(b) Let $B_0=0$ and $B_n=\frac{n(n+1)^{1-\alpha}}{2-\alpha}$. We show
$B_n-B_{n-1}\ge n(n+1)^{-\alpha}$ for every $n\ge1$, which after summing over
$n\le N$ yields $B_N\ge\sum_{n=1}^{N}n(n+1)^{-\alpha}$.

By Bernoulli's inequality, $(1+\frac1n)^\alpha\le1+\frac\alpha n$, so
\begin{align*}
  (n+1)^\alpha\big(B_n-B_{n-1}\big)
  &=\frac{n}{2-\alpha}\Big[n+1-(n-1)\Big(1+\frac1n\Big)^\alpha\Big] \\
  &\ge\frac{n}{2-\alpha}\Big[n+1-(n-1)\Big(1+\frac\alpha n\Big)\Big] \\
  &=\frac{n}{2-\alpha}\Big[2-\alpha+\frac{\alpha}{n}\Big] \\
  &\ge n .
\end{align*}
Dividing by $(n+1)^\alpha$ gives $B_n-B_{n-1}\ge n(n+1)^{-\alpha}$, and summing
over $n=1,\dots,N$ completes the proof.
\end{proof}

\begin{lemma}\label{lem:probit-taylor}
For $q(t)=\Phi^{-1}\big(\tfrac{1+t}{2}\big)$ we have
\begin{equation}\label{eq:probit-series}
  q(t)=\sum_{j\ge1}a_j t^j, \qquad t\in(-1,1),
\end{equation}
with $a_j\ge0$. Moreover,
\begin{equation}\label{eq:probit-bound}
  a_j\le q(\varrho)\,\varrho^{-j}
  \qquad\text{for all } \varrho\in(0,1),
\end{equation}
and the coefficients are given by
\begin{equation}\label{eq:probit-coefficients2}
  a_1=\sqrt{\pi/2},\qquad
  a_{j+1}=\frac{\sqrt{\pi/2}}{j+1}\,[t^j]\exp\Big(\tfrac12\Big(\sum_{i=1}^j a_i t^i\Big)^2\Big).
\end{equation}
\end{lemma}

\begin{proof}
Differentiating the identity $\Phi(q(t))=\tfrac{1+t}{2}$ gives
\begin{equation}\label{eq:qprime}
  q'(t)=\frac{1}{2\varphi(q(t))}=\sqrt{\pi/2}\,e^{q(t)^2/2},
\end{equation}
where $\varphi=\Phi'$ is the standard normal density. Differentiating once more,
\begin{equation}\label{eq:qsecond}
  q''(t)=q(t)\,q'(t)^2,
\end{equation}
and \eqref{eq:qprime} yields the recursion \eqref{eq:probit-coefficients2} at once.

Now consider polynomials in $(q,q')$ with nonnegative coefficients. By the chain rule
and \eqref{eq:qsecond},
\[
  \frac{d}{dt}P(q,q')=\partial_qP\cdot q'+\partial_{q'}P\cdot q\,q'^2,
\]
so the derivative of such a polynomial is again of the same form. Since $q$ and $q'$
are themselves of this form, induction on the order of differentiation shows that
every $q^{(n)}$ is a polynomial in $(q,q')$ with nonnegative coefficients. At $t=0$
we have $q(0)=0$ and $q'(0)=\sqrt{\pi/2}>0$, hence $q^{(n)}(0)\ge0$ and therefore
\[
  a_j=\frac{q^{(j)}(0)}{j!}\ge0.
\]
Also, $q(-t)=\Phi^{-1}\!\big(\tfrac{1-t}{2}\big)=-\Phi^{-1}\!\big(\tfrac{1+t}{2}\big)=-q(t)$,
so $q$ is odd; comparing coefficients in $q(-t)=-q(t)$ gives $a_j=0$ for even $j$.

It remains to identify the series with $q$ on all of $(-1,1)$. Since $\Phi$ is
real-analytic with $\Phi'>0$, the analytic inverse function theorem shows that $q$ is
real-analytic on $(-1,1)$. Let $\varrho_0$ be the radius of convergence of
\eqref{eq:probit-series} and suppose $\varrho_0<1$. By Pringsheim's theorem, a power
series with nonnegative coefficients and finite radius of convergence $\varrho_0$
cannot be analytically continued to a neighbourhood of $t=\varrho_0$. But $q$ is
real-analytic on $(-1,1)\supset(-\varrho_0,\varrho_0]$, which gives exactly such a
continuation - a contradiction. Hence $\varrho_0\ge1$, and by the identity theorem
the series \eqref{eq:probit-series} coincides with $q$ on $(-1,1)$.

Finally, since $a_j\ge0$, for any $\varrho\in(0,1)$,
\[
  a_j\varrho^{j}\le\sum_{i\ge1}a_i\varrho^{i}=q(\varrho),
\]
which is \eqref{eq:probit-bound}.
\end{proof}

\begin{lemma}\label{lem:probit-estimator}
Let $q(t)=\Phi^{-1}\big(\tfrac{1+t}{2}\big)$ and let $a_j$ be the coefficients of its power
series $q(t)=\sum_{j\ge1}a_jt^j$. Fix $m_0\in[0,1)$ and let
\[
  r_1:=\tfrac{1+m_0}2,\quad r_2:=\tfrac{1+r_1}2,\quad \theta:=\tfrac{r_1}{r_2},\quad
  k:=\Big\lceil\tfrac{1-m_0^2}{r_1^2-m_0^2}\Big\rceil,\quad p_j:=(1-\theta)\theta^{j-1}\ \ (j\ge1).
\]
Fix $y\ne z$ with mean answer $\nu$, $|\nu|\le m_0$. Draw $J$ with $\Pr(J=j)=p_j$, make $kJ$ fresh
calls at $(y,z)$, split them into $J$ batches of $k$ calls, let $B_i$ be the average
answer in batch $i$, and set
\[
  T:=\frac{a_J}{p_J}\prod_{i=1}^JB_i .
\]
Then
\begin{enumerate}
\item[(i)] $\E T=q(\nu)$;
\item[(ii)] $\E T^2\le V:=\dfrac{q(r_2)^2\theta^2}{(1-\theta)^2}$;
\item[(iii)] the number of calls $K=kJ$ satisfies $\E e^{uK}=B_0:=1+\theta^{-1/2}$ for
$u:=\dfrac{\log(1/\theta)}{2k}$.
\end{enumerate}
All these constants depend only on $m_0$.
\end{lemma}

\begin{proof}
Let $O_{i,1},\dots,O_{i,k}\in\{-1,+1\}$ be the answers in batch $i$, so that
$B_i=\frac1k\sum_{l\le k}O_{i,l}$. All answers are made at the same pair $(y,z)$ with
fresh noise, hence they are i.i.d.\ with mean $\nu$. They are also independent of $J$,
which is drawn before the calls. Consequently the $B_i$ are independent of each other
and of $J$, and $\E B_i=\nu$.

Since $O_{i,l}^2\equiv1$, we have $\operatorname{Var}(O_{i,l})=1-\nu^2$, so
\[
  \E B_i^2=\operatorname{Var}(B_i)+(\E B_i)^2=\frac{1-\nu^2}k+\nu^2
  =\frac1k+\Big(1-\frac1k\Big)\nu^2 .
\]
The right-hand side is nondecreasing in $\nu^2$. Using $|\nu|\le m_0$
and the choice of $k$,
\begin{equation}\label{eq:batch-second-moment}
  \E B_i^2\le m_0^2+\frac{1-m_0^2}k\le r_1^2 .
\end{equation}
The batch size $k$ is chosen exactly so that this holds: each factor of $T$ then has
second moment at most $r_1^2<1$. By Jensen's inequality, \eqref{eq:batch-second-moment}
also gives $\E|B_i|\le(\E B_i^2)^{1/2}\le r_1$.

\emph{(i)} Given $J=j$, the factors $B_1,\dots,B_j$ are independent with mean $\nu$, so
\[
  \E[T\mid J=j]=\frac{a_j}{p_j}\prod_{i=1}^j\E B_i=\frac{a_j}{p_j}\,\nu^j .
\]
Averaging over $J$,
\[
  \E T=\sum_{j\ge1}p_j\,\frac{a_j}{p_j}\,\nu^j=\sum_{j\ge1}a_j\nu^j=q(\nu)
\]
by \cref{lem:probit-taylor}. Exchanging sum and expectation is justified by absolute
convergence: since $a_j\ge0$ and $\E|B_i|\le r_1$,
\[
  \sum_{j\ge1}p_j\,\E\big[|T|\,\big|\,J=j\big]\le\sum_{j\ge1}a_jr_1^j<\infty ,
\]
because $r_1<1$ lies inside the interval of convergence.

\emph{(ii)} Given $J=j$, independence and \eqref{eq:batch-second-moment} give
\[
  \E[T^2\mid J=j]=\frac{a_j^2}{p_j^2}\prod_{i=1}^j\E B_i^2\le\frac{a_j^2}{p_j^2}\,r_1^{2j}.
\]
Averaging over $J$ and using $a_j\le q(r_2)\,r_2^{-j}$ from \cref{lem:probit-taylor},
$r_1/r_2=\theta$ and $p_j=(1-\theta)\theta^{j-1}$,
\[
  \E T^2\le\sum_{j\ge1}\frac{a_j^2r_1^{2j}}{p_j}
  \le\frac{q(r_2)^2}{1-\theta}\sum_{j\ge1}\frac{\theta^{2j}}{\theta^{j-1}}
  =\frac{q(r_2)^2}{1-\theta}\sum_{j\ge1}\theta^{j+1}
  =\frac{q(r_2)^2}{1-\theta}\cdot\frac{\theta^2}{1-\theta}=V .
\]
The weights $p_j$ decay more slowly than $a_j^2r_1^{2j}\lesssim\theta^{2j}$, which is
why the sum is finite — this is precisely the role of the intermediate radius
$r_2\in(r_1,1)$.

\emph{(iii)} By the choice of $u$, $e^{uk}=\theta^{-1/2}$, hence
$\theta e^{uk}=\theta^{1/2}<1$ and the geometric series converges:
\[
  \E e^{uK}=\E e^{ukJ}=\sum_{j\ge1}(1-\theta)\theta^{j-1}e^{ukj}
  =\frac{(1-\theta)\,e^{uk}}{1-\theta e^{uk}}
  =\frac{(1-\theta)\,\theta^{-1/2}}{1-\theta^{1/2}}
  =\big(1+\theta^{1/2}\big)\theta^{-1/2}=1+\theta^{-1/2}.
\]

All quantities $r_1,r_2,\theta,k,p_j,V,u,B_0$ depend only on
$m_0$, and the coefficients $a_j$ are absolute constants.
\end{proof}

\begin{lemma}[{\cite[Theorem~2]{fokkema2024online}}]\label{lem:value}
Let $Q\subset\R^d$ be a convex body. There are universal constants $C_0,p>0$ such that $\mathsf{BCO}$ (Algorithm 1 of \citet{fokkema2024online}) has the following property. For any
convex $\ell:Q\to[0,1]$, if the observation in round $t$ is $y_t=\ell(x_t)+\zeta_t$,
where $x_t\in Q$ is the point played by the algorithm and, conditionally on 
$x_1,y_1,\dots,x_{t-1},y_{t-1},x_t$, the noise satisfies $\E[\zeta_t\mid\cdot]=0$ and
$\E[e^{\zeta_t^2}\mid\cdot]\le2$, then with probability at least $1-\delta$,
\[
  \sum_{t=1}^n\big(\ell(x_t)-\min_Q\ell\big)
  \le C_0\,M_Q\,d^{2}\sqrt n\,\log^p\!\big(edn/\delta\big),
\]
where $M_Q\in[d^{-1/2},d^{-1/4}]$.
\end{lemma}

\section{EXPERIMENTAL DETAILS}
\label{app:experiments}
All experiments run on a single CPU core.

\paragraph{Non-identifiability (\cref{fig:experiments}, upper left).}
Instances of \cref{thm:nonident} with $d=2$, $R=M=1$, $\mu=0.05$, i.e.\
$f_\sigma(x,\xi)=(1+\tfrac34\sigma\xi)\,q(x,\xi)$ with $\xi$ Rademacher; $30$ seeds; $N=2^{16}$.
The comparison method is \cref{alg:ardg} ($m=1$) fed with
$\sign(f(Y,\xi)-f(Z,\xi))$ for a shared $\xi$ (ties broken by a fair coin), with
$h_k=\min\{0.25,0.5(k+1)^{-1/4}\}$, $\gamma_k=0.3/\sqrt{k+1}$, and uniform averaging.
The transcript does not depend on $\sigma$, so one trajectory per seed is evaluated on both
$f_{+1}$ and $f_{-1}$. Shaded regions: $\pm2$ standard errors for $f_{\pm1}$ and a paired
bootstrap $95\%$ interval for the maximum. The value baseline uses the same seeds
and replaces the bit by the real difference $f_{+1}(Y,\xi)-f_{+1}(Z,\xi)$ divided by $0.5$.

\paragraph{Probit (\cref{fig:experiments}, upper right).}
$d=2$, $Q=\B^2$, $x^*=0.4\,\mathbf 1/\sqrt2$, $F(x)=\frac M{\sqrt2}\|x-x^*\|_1$ with
$M=\frac12$, $\sigma=1$, $R=r=1$, so that $k=3$, $\theta\approx0.818$, $V\approx40.8$.
For each budget $b\in\{2^{10},\dots,2^{19}\}$
we run configuration $\mathrm c$ of \eqref{eq:probit-configurations} with untuned
constants, $G=d\sigma K_0$ with $K_0=\max\{1,M/\sigma,\sqrt V\}$, $\Gamma=M(1+R/r)$,
$h=\min\{r/2,\,RG/(\Gamma\sqrt n)\}$, $\eta=R/(G\sqrt n)$, projection onto $(1-h/r)Q$, and
$n=\lfloor b(1-\theta)/(1.3k)\rfloor$ updates, so that the expected cost is
$b/1.3$; we report the error of the averaged iterate over $10$ seeds. The series for
$\Phi^{-1}((1+t)/2)$ is truncated at order $100$, and the cap of
\cref{alg:probit} is not enforced in this experiment. The strongly convex configuration
has constant $3G^2/\mu$, which puts its asymptotic regime beyond our compute budget, so
we do not report it.
\paragraph{Transfer functions (\cref{fig:experiments}, bottom).}
$d=5$, $Q=\B^5$, $x^*=0.4\,\mathbf 1/\sqrt5$,
$f(x)=\frac\mu2\|x-x^*\|^2+5^{-1/2}\|x-x^*\|_1$ with $\mu=0$ (convex) and $\mu=1$
(strongly convex); $20$ seeds; $N\in\{2^6,\dots,2^{14}\}$; $m=1$;
$h_k=\min\{0.05,\,0.5(k+1)^{-\beta}\}$ with $\beta$ as in
\cref{thm:ardg-convex,thm:ardg-sc} for the respective $q$; convex steps
$\gamma_k=2/\sqrt{k+1}$ with uniform averaging; strongly convex steps
$\gamma_k=\min\{2d/(\lambda\mu h_k(k+1)),\,0.25\}$ with weights $k+1$.
With $t=f(Y)-f(Z)$, the oracle returns $\sign(t+\zeta)$ for symmetric noise $\zeta$:
Gaussian, $\zeta\sim\mathcal N(0,2\sigma^2)$, $\sigma=0.15$ ($\lambda\approx3.76$, $q=2$);
Laplace with scale $a=0.12$ ($\lambda\approx8.33$, $q=1$);
generalized Gaussian with density $\propto e^{-(|z|/a)^\nu}$, $\nu=1.5$, $a=0.12$
($\lambda\approx9.23$, $q=\nu$);
two-sided gamma, $\zeta=\xi_1-\xi_2$ with $\xi_i$ i.i.d.\ symmetric, $|\xi_i|\sim\mathrm{Gamma}(\alpha,\beta)$,
$\alpha=\frac34$, $\beta=12$ ($\lambda\approx10.0$, $q=\frac12$).
Logistic links $\Pr(O=1)=\sigma((t+b)/s)$: Bradley--Terry with $s=0.08$, $b=0$
($\lambda=6.25$, $q=2$); position-biased with $s=0.05$, $b=0.03$ ($\lambda\approx9.15$, $q=1$).
Slopes are least-squares fits of $\log$ error on $\log N$ over the second half of the
checkpoints. Reference lines are anchored at the first point of the Gaussian curve, so only
their slopes are meaningful.

\end{document}